\documentclass[brochure,12pt]{bourbaki}
\usepackage[matrix,arrow]{xy}
\usepackage{amssymb,amsfonts,amsmath,footnote}
\usepackage{mathrsfs}  
\usepackage[utf8]{inputenc}
\usepackage[francais]{babel}
\usepackage{color}
\usepackage{graphicx}
\usepackage{url}

\newcommand{\CC}{\mathbf{C}}
\newcommand{\QQ}{\mathbf{Q}}
\newcommand{\Rep}{\mathbf{Rep}}
\newcommand{\tr}{\mathrm{Tr}}
\newcommand{\GG}{\mathbb{G}}
\newcommand{\FF}{\mathbf{F}}
\newcommand{\alg}{\overline{\mathbf{Q}}}
\newcommand{\cF}{\mathscr{F}}
\newcommand{\cG}{\mathscr{G}}
\newcommand{\IA}{\mathbb{A}}
\newcommand{\cL}{\mathscr{L}}
\newcommand{\pr}{\mathrm{pr}}
\newcommand{\Gal}{\mathrm{Gal}}
\newcommand{\GL}{\mathrm{GL}}
\DeclareMathOperator{\Spec}{Spec}
\newcommand{\Frob}{\mathrm{Frob}}
\newcommand{\arith}{\mathrm{arith}}
\newcommand{\geom}{\mathrm{geom}}
\newcommand{\rg}{\mathrm{rg}}
\newcommand{\Sw}{\mathrm{Sw}}
\newcommand{\PP}{\mathbb{P}}
\newcommand{\Kl}{\mathrm{Kl}}
\newcommand{\Fib}{\mathrm{Fib}}
\newcommand{\ZZ}{\mathbf{Z}}
\newcommand{\cH}{\mathscr{H}}
\newcommand{\Perv}{\mathbf{Perv}}

\newcommand{\DD}{\mathbb{D}}

\newcommand{\bP}{\mathbf{P}}
\newcommand{\dr}{\mathrm{dr}}
\newcommand{\Bon}{\mathrm{Bon}}
\newcommand{\FT}{\mathrm{FT}}

\newcommand{\pct}{\mathrm{pct}}
\newcommand{\npct}{\mathrm{npct}}

\newcommand{\Mauv}{\mathrm{Mauv}}
\newcommand{\gr}{\mathrm{gr}}
\newcommand{\cO}{\mathcal{O}}
\newcommand{\SL}{\mathrm{SL}}

\newcommand{\RR}{\mathbf{R}}
\newcommand{\Sp}{\mathrm{Sp}}
\newcommand{\sep}{\textup{sep}}
\newcommand{\mongeo}{G_{\geom,\hspace{.3mm} M}}
\newcommand{\monarith}{G_{\arith,\hspace{.3mm} M}}

\usepackage{bbm}
\usepackage{tensor}
\usepackage{enumitem}

\let\originalleft\left
\let\originalright\right
\def\left#1{\mathopen{}\originalleft#1}
\def\right#1{\originalright#1\mathclose{}}

\newsavebox{\smlmat}
\savebox{\smlmat}{$\left(\begin{smallmatrix}7&-2\\4&-1\end{smallmatrix}\right)$}

\usepackage{accents}
\newcommand{\ubar}[1]{\underaccent{\bar}{#1}}

\let\leq\leqslant
\let\geq\geqslant

\date{Janvier 2018}
\bbkannee{70\`eme ann\'ee, 2017-2018}
\bbknumero{1141}
\title{\'Equirépartition de sommes exponentielles}
\subtitle{Travaux de Katz}

\author{Javier FRES\'AN}
\address{CMLS, \'Ecole Polytechnique\\ 
F-91128 Palaiseau Cedex}
\email{javier.fresan@polytechnique.edu}

\begin{document}

\maketitle

\section{Introduction}

Cent ans se sont écoulés depuis la parution de l'article de Weyl \cite{weyl-equi} sans que les résul\-tats d'équirépartition  en théorie des nombres cessent de faire florès, révélant des liens profonds avec la géométrie algébrique et la théorie des représentations. L'interaction de ces domaines est particulièrement riche dans les travaux de Katz dont il sera question ici. Ce sont des théorèmes, anciens et récents, d'équirépartition  de sommes exponentielles sur les corps finis, le plus souvent à caractéristique fixée. Les sommes concernées s'obtiennent par transformation de Fourier, relative à un caractère, de la fonction trace d'un faisceau $\ell$-adique sur un groupe algébrique commutatif, et il s'agit de comprendre leur répartition lorsque le faisceau est fixe mais que l'on fait varier le~caractère.  

\subsection{L'exemple des sommes de Gauss et des sommes de Kloosterman}

Soient $p$ un nombre premier, $q$ une puissance de $p$ et $\FF_q$ un corps fini à $q$ éléments. Étant donnés un caractère additif non trivial $\psi \colon \FF_q \to \CC^\times$ et un caractère multiplicatif~$\chi \colon \FF_q^\times \to \CC^\times$, on définit la \textit{somme de Gauss} $g(\psi, \chi)$ comme l'entier algébrique
\begin{displaymath}
g(\psi, \chi)=\sum_{x \in \FF_q^\times} \psi(x) \chi(x). 
\end{displaymath} Par exemple, si $q=p$ et si l'on prend pour $\psi$ le caractère $x \mapsto \exp(2\pi ix \slash p)$ et pour $\chi$ le symbole de Legendre, il s'agit de la somme considérée par Gauss dans sa quatrième preuve de la loi de réciprocité quadratique~\cite{gauss}. 

Si $\chi$ est trivial, la somme de Gauss vaut $-1$; sinon, sa valeur absolue est égale à $\sqrt{q}$. Choisissons, pour chaque $p$, un caractère non trivial~$\psi_p$ de~$\FF_p$ et notons $\psi_q$ le caractère de $\FF_q$ obtenu par composition avec la trace. Gardant $\psi_q$ fixe et faisant varier $\chi$ parmi les caractères multiplicatifs non triviaux, on obtient~$q-2$ {points
\begin{displaymath}
\theta_{q, \chi}=\frac{g(\psi_q, \chi)}{\sqrt{q}}
\end{displaymath} dans le cercle unité $S^1$. Comment ces points se répartissent-ils quand $q$ tend vers l'infini ?

Soient $(X, \mu)$ un espace topologique compact muni d'une mesure de probabilité~$\mu$ et~$(S_N)$ une suite\footnote{Bien que l'on utilise $N$ comme paramètre, dans les énoncés qui suivent la suite n'est pas forcément indexée par des entiers mais par des caractères, des points d'une variété algébrique, etc.} d'ensembles finis non vides avec des applications~\hbox{$\theta_N \colon S_N \to X$.} Rappelons que les $(S_N, \theta_N)$ sont dits \textit{équirépartis  selon~$\mu$} si la suite de mesures~$|S_N|^{-1}\sum_{x \in S_N} \delta_{\theta_N(x)}$ converge vaguement vers $\mu$ lorsque $N$ tend vers l'infini, c'est-à-dire si, pour toute fonction continue $f \colon X \to \CC$, on a l'égalité 
\begin{equation*}\label{eqn:equidis} 
\int_X f(x)\mu(x)=\lim_{N \to \infty} \frac{1}{|S_N|} \sum_{x \in S_N} f(\theta_N(x)). 
\end{equation*} Il suffit en fait de la vérifier pour une classe de fonctions test dont les combinaisons linéaires finies sont denses dans l'espace $\mathscr{C}(X)$ des fonctions continues à valeurs complexes muni de la topologie de la convergence uniforme.

Dans le cas qui nous occupe, Katz remarqua dans \cite[\S1.3.3]{Katz80} que la majoration des sommes de Kloosterman obtenue par Deligne comme conséquence de ses travaux sur la conjecture de Weil entraînait le résultat d'équirépartition  suivant:

\begin{theo}[Deligne]\label{theo-gauss} Lorsque $q$ tend vers l'infini, les points $\{ \theta_{q, \chi}\}_{\chi \neq 1}$ s'équirépartissent selon la mesure de Haar normalisée sur le cercle unité. Autrement dit, pour toute fonction continue $f \colon S^1 \to \CC$ on a l'égalité
\begin{equation}\label{eqn:gauss}
\frac{1}{2\pi} \int_0^{2\pi} f(e^{i\theta})\,d\theta=\lim_{q \to \infty} \frac{1}{q-2} \sum_{\chi \neq 1} f(\theta_{q, \chi}).
\end{equation}
\end{theo}

Comme les polynômes de Laurent sont denses dans $\mathscr{C}(S^1)$, il suffit de considérer les fonctions $f(z)=z^n$ avec $n$ entier. Le cas $n=0$ est évident. Le membre gauche de \eqref{eqn:gauss} étant nul pour $n \neq 0$, il faut démontrer que la suite
\begin{displaymath}\label{eqn:moments}
\frac{1}{q-2}\sum_{\chi \neq 1} f(\theta_{q, \chi})= \frac{1}{q^{\frac{n}{2}}(q-2)} \sum_{\chi \neq 1} g(\psi_q, \chi)^n
\end{displaymath} converge vers zéro lorsque $q$ tend vers l'infini (c'est le \textit{critère d'équirépartition  de Weyl}). Grâce à la relation $g(\psi_q, \chi)^{-1}=g(\bar{\psi}_q, \bar{\chi})q^{-1}$, on peut se ramener à $n \geq 1$, auquel cas les puissances des sommes de Gauss sont égales à \vspace{1mm}
\begin{align}\label{eqn:fourier-gauus}
g(\psi_q, \chi)^n&=\sum_{x_1, \dots, x_n \in \FF_q^\times} \psi_q(x_1+\cdots+x_n)\chi(x_1\cdots x_n) \nonumber \\
&=\sum_{a \in \FF_q^\times} \chi(a) \sum_{\substack{x_1, \dots, x_n \in \FF_q^\times \\ x_1\cdots x_n=a}} \psi_q(x_1+\cdots+x_n).
\end{align} 

Il s'ensuit que $\chi \mapsto g(\psi_q, \chi)^n$ est la transformée de Fourier, au sens du groupe abélien fini $\FF_q^\times$, de la fonction qui à un élément $a \in \FF_q^\times$ associe la \textit{somme de Kloosterman} 
\begin{equation*}\label{eqn:somme-kloos}
\Kl_n(a, q)=\sum_{\substack{x_1, \dots, x_n \in \FF_q^\times \\ x_1\cdots x_n=a}} \psi_q(x_1+\ldots+x_n).
\end{equation*} Pour $n=2$, ces sommes font un caméo dans l'article posthume de Poincaré sur les formes modulaires, où il \og se borne à constater \fg{} qu'elles ne sont pas nulles en général\footnote{Les sommes de Poincaré portent sur $(\ZZ/q\ZZ)^\times$ et peuvent être nulles. Celles que l'on considère ici ne le sont \textit{jamais}, car $\Kl_2(a, q)$ appartient au sous\nobreakdash-anneau de $\CC$ engendré par une racine primitive~$p$\nobreakdash-ème de l'unité $\zeta_p$ et l'on a~$\Kl_2(a, q) \equiv -1$ modulo l'idéal premier $1-\zeta_p$.}~\cite[p.\,148]{poincare}. Kloosterman les introduisit de manière indépendante en 1926, en raffinant la~méthode du cercle pour étudier l'asymptotique du nombre de représentations d'un entier par une forme quadratique définie positive en quatre variables~\cite{kloosterman}. Un point clé de son travail est la majoration $|\Kl_2(a, q)| < 2q^{3/4}$, qu'il obtint en calculant le quatrième moment
\begin{equation}\label{eqn:4ememoment}
\sum_{a \in \FF_q^\times} \Kl_2(a, q)^4=2q^3-3q^2-3q-1.
\end{equation} Quelques années plus tard, Salié~\cite{salie} et Davenport~\cite{davenport-kloos} purent améliorer l'exposant de~$3/4$ à $2/3$ en estimant le sixième moment. Puis en 1934 Hasse observa, en comparant la somme de Kloosterman au nombre de solutions de l'équation $y^q-y=x+ax^{-1}$, que la borne optimale $2 \sqrt{q}$ découlait de l'hypothèse de Riemann pour les courbes sur les corps finis \cite{hasse}; avec la preuve de Weil entre 1940 et 1948, elle fut enfin établie \cite{weil-exp}. 

Trouver la majoration optimale pour les sommes de Kloosterman en plusieurs variables est une tâche significativement plus compliquée qui requiert l'analogue de l'hypothèse de Riemann pour la cohomologie à coefficients dans un faisceau $\ell$-adique. En la démontrant dans l'article \cite{weil2}, auquel on se référera comme \og Weil II \fg{} par la suite, Deligne ouvrit la voie à de nombreuses applications à l'étude des sommes exponentielles que nous sommes encore loin d'avoir épuisées. Lui-même exposa le principe de la méthode dans~\cite{sommes-trigo}, où il montre comment en déduire l'estimée
\begin{displaymath}
|\Kl_n(a, q)| \leq n q^{\frac{n-1}{2}}
\end{displaymath} pour n'importe quels $n \geq 2$ et $a \in \FF_q^\times$. C'est ce qu'il fallait pour conclure la preuve. 

\noindent {\sc Fin de la démonstration}  --- En sommant la formule \eqref{eqn:fourier-gauus} sur les caractères multiplicatifs non triviaux, il vient
\begin{align*}
\sum_{\chi \neq 1} g(\psi_q, \chi)^n&=-g(\psi, 1)^n+\sum_{a \in \FF_q^\times} \Kl_n(a, q) \sum_\chi \chi(a)=(-1)^{n+1}+(q-1) \Kl_n(1, q)
\end{align*} par orthogonalité des caractères. 
D'après la majoration de Deligne, nous avons donc
\begin{equation*}\label{eqn:majgauss}
\left|\frac{1}{q^{\frac{n}{2}}(q-2)} \sum_{\chi \neq 1} g(\psi_q, \chi)^n\right| \leq \frac{2n+1}{\sqrt{q}}
\end{equation*} pour tout $q>2$ et la limite du membre gauche est bien zéro lorsque $q \to \infty$.
\qed

\begin{rema} Le théorème s'applique tant à la situation où $p$ est fixe et l'on fait tendre $q$ vers l'infini parmi les puissances de $p$ qu'à la situation où $p$ varie aussi, ce qui est possible car la constante dans la majoration des moments (en l'occurrence~$2n+1$) est indépendante de $p$. En théorie analytique des nombres, on parle souvent d'équirépartition \textit{verticale} ou \textit{horizontale} pour distinguer ces deux cas.
\end{rema}

Revenons maintenant aux sommes de Kloosterman. Comme appliquer la conjugaison complexe revient à échanger $x_i$ et $-x_i$ dans l'expression de $\Kl_2(a, q)$, ce sont des nombre réels (il en va de même pour tout $n$ pair). Au vu de la borne de Weil, pour chaque~$a \in \FF_q^\times$, il existe un unique angle $\theta_{q, a} \in [0, \pi]$ tel~que 
\begin{displaymath}
\Kl_2(a, q)=2\sqrt{q}\cos \theta_{q, a}.
\end{displaymath} Comment ces $q-1$ angles varient-ils avec $q$? En s'appuyant sur l'interprétation de la somme~$\Kl_2(a, q)$ comme la trace de Frobenius en $a \in \mathbb{G}_m(\FF_q)$ d'un système local~$\ell$\nobreakdash-adique sur $\GG_m$ et sur un théorème de Deligne affirmant que la répartition de telles traces est gouvernée par le groupe de monodromie, Katz démontra dans \cite{katz-kloos} que les angles~$\{\theta_{q, a}\}_{a \in \FF_q^\times}$ se répartissent comme les classes de conjugaison de matrices aléatoires dans le groupe spécial unitaire~$\mathrm{SU}(2)$. Plus précisément, si l'on identifie l'intervalle $[0, \pi]$ à l'espace de ces classes par l'application qui envoie $\theta$ sur la classe de conjugaison de~$\left(\begin{smallmatrix} e^{i \theta} & 0 \\ 0 & e^{-i\theta} \end{smallmatrix}\right)$, l'image directe de la mesure de Haar normalisée sur~$\mathrm{SU}(2)$ par la projection canonique est la mesure $(2 \slash\pi) \sin^2\theta \,d\theta$ sur~$[0, \pi]$. Elle porte en théorie des nombres le nom de \textit{mesure de~Sato--Tate}, d'après leur célèbre conjecture sur la répartition du terme d'erreur dans l'approximation par~$p+1$ du nombre de points de la réduction modulo $p$ d'une courbe elliptique sur $\QQ$ sans multiplication complexe.  

\begin{theo}[Katz]\label{thm:equi-kloos} Lorsque $q$ tend vers l'infini, les angles $\{ \theta_{q, a}\}_{a \in \FF_q^\times}$ s'équirépartissent selon la mesure de Sato--Tate, c'est-à-dire pour toute fonction continue $f \colon [0, \pi] \to \CC$ on a l'égalité 
\begin{displaymath}
\frac{2}{\pi} \int_0^\pi f(\theta) \sin^2 \theta \, d\theta=\lim_{q \to \infty} \frac{1}{q-1} \sum_{a \in \FF_q^\times} f(\theta_{q, a}). 
\end{displaymath}
\end{theo}

À nouveau, ce résultat est valable pour n'importe quelle suite de corps finis de cardinaux croissants. On conjecture également que, pour un entier \textit{fixe} $a$, les angles~$\{\theta_{p, a}\}_{p \leq N}$ s'équirépartissent selon la mesure de Sato--Tate lorsque $p$ tend vers l'infini parmi les nombres premiers ne divisant pas $a$, mais même des conséquences faibles de cet énoncé paraissent hors de portée à l'heure actuelle\footnote{Mentionnons, à titre d'exemple, la \textit{conjecture du changement de signe}: puisque $\sin^2 \theta$ est symétrique par rapport à $\theta=\pi/2$ et que les $\Kl_2(a, p)$ sont tous non nuls, il devrait y avoir asymptotiquement autant de premiers $p$ pour lesquels la somme de Kloosterman est positive que négative. On n'en sait rien! Dans \cite{fouvry-michel}, Fouvry et Michel démontrent des résultats dans cette direction quand la suite des $p$ est remplacée par une suite de nombres \og presque premiers \fg{}, c'est-à-dire non premiers, sans facteur carré et ayant un nombre de facteurs premiers borné uniformément.}. 

\subsection{Les sommes d'Evans et de Rudnick} 

Autour de 2003, Evans et Rudnick trouvèrent d'autres exemples de sommes exponentielles dépendant d'un caractère multiplicatif qui, d'après leurs expériences numériques, semblaient s'équirépartir selon la mesure de Sato--Tate. D'un côté, Evans étudia les~sommes
\begin{displaymath}
S(\chi)=-\frac{1}{\sqrt{q}} \sum_{x \in \FF_q^\times} \chi(x) \psi_q\left(x-x^{-1}\right), 
\end{displaymath} qui sont des nombres réels de valeur absolue au plus~$2$ par l'hypothèse de Riemann sur les courbes. Ils s'écrivent donc $S(\chi)=2\cos \theta_{q, \chi}$ pour un unique $\theta_{q, \chi} \in [0, \pi].$ Est-il vrai que les angles $\{\theta_{q, \chi} \}_{\chi}$ s'équirépartissent selon la mesure de Sato--Tate quand $q$ tend vers l'infini? Ou, ce qui revient au même, les sommes~$S(\chi)$ s'équirépartissent-elles selon la \textit{mesure du demi-cercle} $(1\slash 2\pi)\sqrt{4-x^2}\, dx$ sur l'intervalle $[-2, 2]$?

À la même époque, Rudnick rencontra des sommes semblables dans ses travaux sur le chaos quantique. Décrivons brièvement le contexte, en nous référant à \cite{hecke-cat}, \cite{rudnick-annals} et \cite{quantumcat} pour plus de détails. Il s'agit de quantifier le système dynamique obtenu en itérant la transformation du tore $T^2=\RR^2 \slash \ZZ^2$ définie par une matrice hyperbolique $A \in \SL(2, \ZZ)$, des applications connues comme \og cat maps \fg{} dans la littérature. 
Soient $N \geq 1$ un entier et $\mathcal{H}_N$ l'espace de Hilbert $L^2(\ZZ/N\ZZ)$. On pense à $N$ comme à l'inverse de la constante de Planck $\hbar$, de sorte que la limite semi-classique $\hbar \to 0$ devienne $N \to \infty$. Après avoir associé à chaque observable classique $f \colon T^2 \to \RR$ de classe $\mathscr{C}^\infty$ un opérateur autoadjoint $\mathrm{Op}_N(f)$ sur $\mathcal{H}_N$ (l'observable quantique), une quantification de $A$ est une suite d'opérateurs unitaires $U_N(A)$ sur $\mathcal{H}_N$ vérifiant la formule d'Egoroff 
$$
U_N(A)^\ast\mathrm{Op}_N(f)U_N(A)=\mathrm{Op}_N(f \circ A)
$$ pour toute observable $f$. Si $\varphi_1, \dots, \varphi_N$ est une base orthonormale de $\mathcal{H}_N$ formée de vecteurs propres pour l'opérateur $U_N(A)$, on sait que la moyenne des espé\-rances~$\langle \mathrm{Op}_N(f) \varphi_i, \varphi_i \rangle$ converge vers $\int_{T^2} f(x) dx$ lorsque $N$ tend vers l'infini. Est-ce encore vrai pour chaque terme pris individuellement?

Dans \cite{hecke-cat}, Kurlberg et Rudnick observent que $A$ admet une quantification qui ne dépend que de la réduction de la matrice modulo $2N$, ce qui leur permet d'introduire une classe d'opérateurs d'origine arithmétique sur $\mathcal{H}_N$ qu'ils appellent opérateurs de Hecke. Le théorème d'\textit{unique ergodicité quantique} est l'énoncé que, si les $\varphi_i$ sont simultanément des fonctions propres pour $U_N(A)$ et les opérateurs de Hecke, on a
\begin{displaymath}
\lim_{N \to \infty} \langle \mathrm{Op}_N(f) \varphi_i, \varphi_i \rangle=\int_{T^2} f(x) \,dx
\end{displaymath} pour tout $i$. Comment les espérances $\langle \mathrm{Op}_N(f) \varphi_i, \varphi_i \rangle$ fluctuent-elles autour de la limite? Puisque l'on s'attend à ce que le terme d'erreur soit génériquement d'ordre  $1\slash \sqrt{N}$, il est naturel de considérer les fluctuations définies par 
\begin{displaymath}
F_i^{(N)}=\sqrt{N} \left(\langle  \mathrm{Op}_N(f) \varphi_i, \varphi_i\rangle- \int_{T^2} f(x) \, dx\right).
\end{displaymath} 

Prenons pour $N$ un nombre premier $p \geq 3$ tel que la matrice $A$ soit diagonalisable modulo~$p$, disons par une matrice de passage $P$. Dans ce cas, chaque caractère multiplicatif $\chi$ de~$\FF_p$ donne lieu à une fonction propre $\varphi_\chi=U_N(P) \chi$ pour $U_N(A)$ et les opérateurs de Hecke. À un facteur près, les fluctuations des $\varphi_\chi$ sont les sommes
\begin{displaymath}
E(\chi)=-\frac{1}{\sqrt{p}} \sum_{x \in \FF_p \setminus \{0, 1\}} \chi(x) \psi_p\left((x+1)\slash(x-1)\right).
\end{displaymath} Si $\chi$ n'est pas trivial, $E(\chi)$ est encore un nombre réel compris entre $-2$ et~$2$. Des expériences numériques (figure \ref{fig1}) suggèrent que les $\{E(\chi) \}_{\chi \neq 1}$ s'équirépartissent selon la mesure du demi-cercle lorsque~$p$ tend vers l'infini. Peut-on le démontrer? 

\begin{figure}
\begin{center}
  \includegraphics[width=0.75\textwidth]{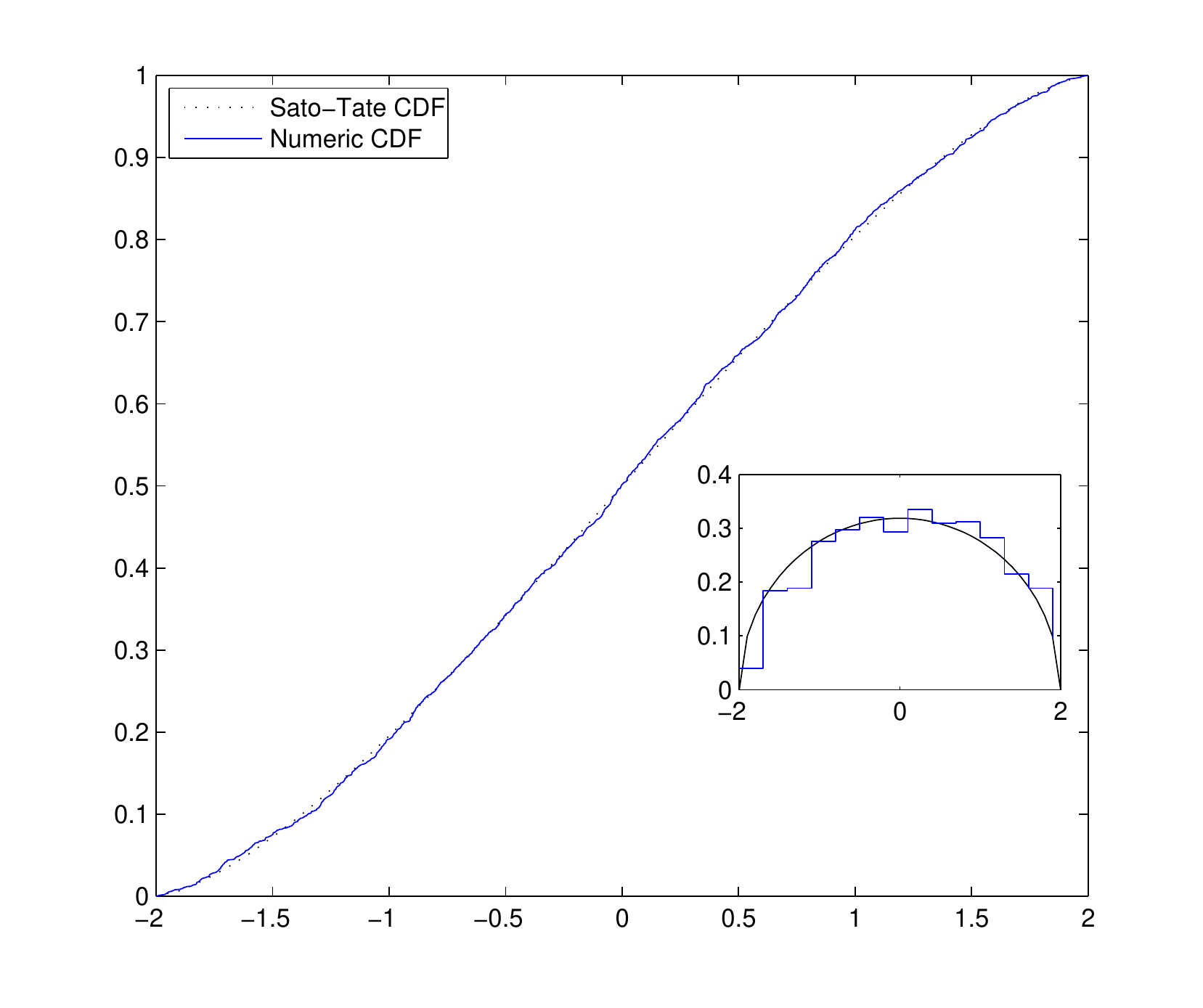}
  \end{center}
  \caption{Comparaison entre la fonction de répartition des fluctuations pour la matrice $A=\usebox{\smlmat}$ et le premier $p=1997$ et celle d'une variable aléatoire pour la mesure du demi-cercle, d'après Rosenzweig \cite{irmn-ros}.}
  \label{fig1}
\end{figure}

Dans ce texte, nous expliquerons comment un théorème de Katz \cite{katz} permet entre autres d'apporter une réponse positive aux questions d'Evans et de Rudnick. Leur difficulté vient du fait que, n'étant pas les traces de Frobenius d'un même système local~$\ell$\nobreakdash-adique, ces familles de sommes exponentielles indexées par des caractères multiplicatifs ne rentrent pas dans le cadre du théorème d'équirépartition de Deligne. Il faut notamment trouver un substitut du groupe de monodromie: c'est là que le formalisme tannakien nous viendra en aide. Avant d'énoncer le théorème principal, nous illustrerons la démarche par les sommes de~Gauss\footnote{La recension \cite{review-kowalski} est un excellent point d'entrée aux idées de la preuve. L'exemple des sommes de Gauss est certes très particulier, en ce sens que pour déterminer leur répartition il suffit d'estimer leurs puissances $n$-ièmes et que l'on arrive à identifier celles-ci avec les sommes de Kloosterman, mais en l'interprétant géométriquement on y voit déjà apparaître tous les ingrédients essentiels.}. 

\subsection{La méthode $\ell$-adique}\label{sec:intromethode}

Posons $k=\FF_q$ et fixons une clôture algébrique $\bar k$ de $k$. Soient $\ell$ un nombre premier distinct de $p$ et $\alg_\ell$ une clôture algébrique du corps des nombres $\ell$-adiques. Le choix d'un plongement $\iota \colon \alg_\ell \hookrightarrow \CC$ induit un isomorphisme entre les clôtures algébriques de~$\QQ$ dans~$\alg_\ell$ et dans $\CC$ par le biais duquel on peut voir les caractères à valeurs dans~$\alg_\ell^\times$. Soit donc~$\psi \colon k \to \alg_\ell^\times$ un caractère additif non trivial. Considérons une clôture séparable~$k(T)^{\sep}$ du corps des fractions rationnelles et une solution~\hbox{$U \in k(T)^{\sep}$} de l'équation d'Artin--Schreier~$U-U^q=T$. L'extension~$k(T, U)$ est galoisienne sur $k(T)$ de groupe~$k$, un élément $a \in k$ opérant par $U \mapsto U+a$. Notons~$\cL_\psi$ la représentation 
\begin{displaymath}
\Gal(k(T)^{\sep} \slash k(T)) \longrightarrow \Gal(k(T, U) \slash k(T)) \cong k \stackrel{\psi}{\longrightarrow} \alg_\ell^\times
\end{displaymath} obtenue par composition de la surjection canonique avec le caractère $\psi$. Identifiant~$k(T)$ au corps de fonctions de la droite projective $\PP^1_k$, on voit que $\cL_\psi$ est non ramifiée en dehors de l'infini. Par suite, pour chaque point $x$ de la droite affine~$\IA^1_k$ à valeurs dans une extension finie $E$ de $k$, la classe de conjugaison $\Frob_{E, x}$ dans le groupe~\hbox{$\Gal(k(T)^{\sep} \slash k(T))$} définie par le Frobenius géométrique de $\Gal(\bar{k} \slash E)$ agit sur~$\alg_\ell^\times$ au travers de la représentation $\cL_\psi$. Pour calculer cette action, notons $n$ le degré de $E$ sur~$k$. Si $y \in \bar{k}$ satisfait à l'équation $y-y^q=x$, alors $y^{q^j}=y^{q^{j-1}}-x^{q^{j-1}}$ pour tout entier $j \geq 1$; sommant sur les $j=1, \dots, n$, on trouve que la substitution de Frobenius est donnée par $y \mapsto y^{q^n}=y-\tr_{E \slash k}(x)$. Le Frobenius géométrique étant son inverse, on en déduit après application du caractère $\psi$ que la classe de conjugaison $\Frob_{E, x}$ agit sur l'espace de la représentation $\cL_\psi$ comme la multiplication par~$\psi(\tr_{E \slash k}(x)).$ 

De même, à partir de l'équation de Kummer $U^{1-q}=T$, on peut associer à chaque caractère multiplicatif~$\chi \colon k^\times \to \alg_\ell^\times$ une représentation $\cL_\chi$ de $\Gal(k(T)^{\sep} \slash k(T))$ à valeurs dans $\alg_\ell^\times$, cette fois-ci ramifiée en zéro et l'infini, telle que la classe de conjugaison~$\Frob_{E, x}$ agisse par multiplication par $\chi(\mathrm{N}_{E \slash k}(x))$ pour tout $x \in E^\times$. En termes de ces représentations, la somme de Gauss s'écrit 
\begin{displaymath}
g(\psi, \chi)=\sum_{x \in \FF_q^\times} \tr(\Frob_{E, x} \hspace{.3mm} | \hspace{.3mm} \cL_\psi \otimes \cL_\chi). 
\end{displaymath} 

Dans le langage plus géométrique que l'on adoptera ici, $\cL_\psi$ et $\cL_\chi$ sont des systèmes locaux $\ell$-adiques de rang un sur $\IA^1_k$ et sur $\GG_{m, k}$ respectivement; on notera encore $\cL_\psi$~la restriction à $\GG_{m, k}.$ Pour des raisons pour l'instant mystérieuses, il est préférable de travailler avec l'objet $M=\cL_\psi(1/2)[1]$. Le symbole~$(1/2)$ désigne une \textit{demi\nobreakdash-torsion à la Tate}, dont l'effet est de multiplier les traces de Frobenius par $1/\sqrt{q}$, et le décalage~$[1]$ indique que l'on regarde $\cL_\psi(1/2)$ non pas comme un faisceau mais comme un complexe concentré en degré $-1$, ce qui change le signe des traces. Ainsi~$M$ est un objet de la catégorie \textit{dérivée} $D^b_c(\GG_{m, k}, \alg_\ell)$ de faisceaux $\ell$\nobreakdash-adiques sur~$\GG_{m, k}$. Grâce au décalage, il appartient en fait à une sous-catégorie \textit{abélienne} formée d'objets avec des propriétés remarquables que l'on appelle \textit{faisceaux pervers}, même s'il s'agit plutôt de~complexes. 

Étant donnés une extension finie $E$ de $k$ et un caractère multiplicatif $\chi \colon E^\times \to \alg_\ell^\times$, l'objet $M \otimes \cL_\chi$ est encore un faisceau pervers sur $\GG_{m, E}$ dont les groupes de cohomologie à support compact s'annulent en degré non nul. En degré zéro,  
\begin{displaymath}
\omega_\chi(M)=H^0_c(\GG_{m, \bar{k}}, M \otimes \cL_\chi)
\end{displaymath} est une $\alg_\ell$-droite munie d'une action du groupe $\Gal(\bar{k} \slash E)$, en particulier du Frobenius géométrique $F_E.$ La formule des traces de Grothendieck est dans ce cas l'égalité 
\begin{equation}\label{eqn:gauss-sum-over-E}
\tr(F_E \hspace{.3mm} | \hspace{.3mm} \omega_\chi(M))  
=\sum_{x \in E^\times} \tr(\Frob_{E, x} \hspace{.3mm} | \hspace{.3mm} M \otimes \cL_\chi)=-\frac{1}{\sqrt{|E|}}\sum_{x \in E^\times} \psi(x) \chi(x).
\end{equation} 

La distinction selon que le caractère $\chi$ soit trivial ou non se reflète dans l'action de Frobenius sur $\omega_\chi(M)$, qui est en l'occurrence réduite à la multiplication par le nombre~$\ell$\nobreakdash-adique \eqref{eqn:gauss-sum-over-E}. À travers le plongement $\iota \colon \alg_\ell \hookrightarrow \CC$, on peut parler de sa valeur absolue: elle est égale à $1$ si et seulement si $\chi$ n'est pas trivial. Pour des objets~$M$ plus généraux, on dira qu'un caractère multiplicatif $\chi$ est \textit{bon} pour $M$ si toutes les valeurs propres de $F_E$ agissant sur $\omega_\chi(M)$ sont unitaires. 

\subsection{Convolution et équirépartition }

L'étape suivante consiste à introduire une opération sur les faisceaux relevant le produit des sommes exponentielles; c'est là que la structure de groupe sur $\GG_{m, k}$ intervient. Pour ce faire, on définit sur la catégorie dérivée $D^b_c(\GG_{m, k}, \alg_\ell)$ un produit de \textit{convolution} 
\begin{displaymath}
M \star_! N=Rm_!(\pr_1^\ast M \otimes \pr_2^\star N), 
\end{displaymath} où $\pr_1$ et $\pr_2$ désignent les projections de $\GG_{m, k} \times \GG_{m, k}$ sur les deux facteurs et $m$ la loi de multiplication. Un premier obstacle est que ce produit ne préserve pas les faisceaux pervers, à moins de se restreindre à une sous-catégorie convenable $\bP(\GG_{m, k})$. Un second point technique est que l'objet $M\star_!N$ doit être remplacé par un quotient~\hbox{$M \star_{\textup{int}} N$} que l'on appelle convolution \textit{intermédiaire}.  
Si $M$ et $N$ sont des faisceaux pervers dans~$\bP(\GG_{m, k})$ pour lesquels un caractère $\chi$ est bon, on a un isomorphisme 
$$
\omega_\chi(M \star_{\textup{int}} N) \cong \omega_\chi(M) \otimes \omega_\chi(N) 
$$ compatible à l'action de Frobenius. D'après Gabber et Loeser \cite{gabber-loeser}, avec la convolution intermédiaire pour produit tensoriel, $\bP(\GG_{m, k})$ a la structure d'une catégorie \textit{tannakienne} dans laquelle l'objet neutre est le faisceau gratte-ciel supporté en $1$, le dual $M^\vee$ est le tiré en arrière du dual de Verdier $\DD(M)$ par l'inversion $x \mapsto x^{-1}$ et les objets ont pour dimension leur caractéristique d'Euler. 

Concrètement, on se servira de cette information de la manière suivante. Soit $\langle M \rangle^{\otimes}$ la sous\nobreakdash-catégorie~de $\bP(\GG_{m, k})$ formée des objets qui s'obtiennent à partir de $M$ en itérant les opérations produit tensoriel, dual, somme directe et sous-quotient, autrement dit qui sont des sous-quotients d'une \textit{construction tensorielle} \hbox{$M^{\ubar{r}, \ubar{s}}=\bigoplus M^{\otimes r_i} \otimes (M^\vee)^{\otimes s_i}$;} c'est la plus petite sous-catégorie tannakienne de $\bP(\GG_{m, k})$ contenant l'objet $M$. De plus,~$\omega_\chi$ est un foncteur fibre sur~$\langle M \rangle^{\otimes}$ pour tout bon caractère $\chi$ et les automorphismes linéaires de $\omega_\chi(M)$ agissent naturellement sur tous les $\omega_\chi(M^{\ubar{r}, \ubar{s}})$. Soit $$G_{M, \chi} \subseteq \GL(\omega_\chi(M))$$ le sous-groupe de ceux laissant globalement stables les images par $\omega_\chi$ de n'importe quel sous-quotient de $M^{\ubar{r}, \ubar{s}}$. C'est un sous-groupe algébrique défini sur $\alg_\ell$ dont la catégorie des représentations de dimension finie est, par la théorie tannakienne, équivalente à~$\langle M \rangle^{\otimes}$. On aura également besoin de prendre en compte le faisceau pervers $M_{\bar k}$ sur~$\GG_{m, \bar k}$ déduit de $M$ par extension des scalaires. Comme ses constructions tensorielles ont plus de sous-quotients à respecter, le groupe correspondant est \textit{a priori} plus petit; pour les distinguer, on notera $G_{\arith,\hspace{.3mm} M, \chi}$ celui de $M$ et $G_{\geom,\hspace{.3mm} M, \chi}$ celui de $M_{\bar k}$ et on les appellera arithmétique et géométrique respectivement.

Le Frobenius $F_E$ est un automorphisme de $\omega_\chi(M)$ vérifiant la propriété de stabilité ci-dessus et définit ainsi un élément $\Frob_{E, \chi}$ dans le groupe $G_{\arith,\hspace{.3mm} M, \chi}$. Si l'on change de caractère, on change de groupe, la dépendance en $\chi$ étant la même que celle du groupe fondamental d'un espace topologique en le point base: tous les $G_{\arith,\hspace{.3mm} M, \chi}$ sont isomorphes entre eux par des isomorphismes uniques à des automorphismes intérieurs près. On peut alors fixer un caractère $\chi_0$ et voir les différents $\Frob_{E, \chi}$ comme des classes de conjugaison dans le même groupe $G_{\arith,\hspace{.3mm} M, \chi_0}$, que l'on notera simplement $\monarith$. On ne sait pas si $\Frob_{E, \chi}$ est diagonalisable mais, quitte à prendre sa semisimplification~$\Frob_{E, \chi}^{\textup{ss}}$ au sens de la décomposition de Jordan, on obtient une classe de conjugaison semisimple dans $\monarith(\alg_\ell)$ avec des valeurs propres unitaires, que l'on peut ensuite regarder dans $\monarith(\CC)$ via le plongement $\iota \colon \alg_\ell \hookrightarrow \CC$. 

Les propriétés des objets d'une catégorie tannakienne se traduisent en des propriétés de leurs groupes et \textit{vice versa}. Par exemple, si $M$ est somme directe d'objets simples, alors le groupe $\monarith$ est réductif car il possède une représentation fidèle et complètement réductible. Supposons que c'est le cas et choisissons un sous-groupe compact maximal $K$ du groupe de Lie complexe $\monarith(\CC)$. Chaque élément dans $\Frob_{E, \chi}^{\textup{ss}}$ est alors conjugué à un élément dans $K$, qui est à son tour bien défini à conjugaison près par un élément de $K$. On en déduit une classe de conjugaison $\theta_{E, \chi}$ dans $K$ dont les traces $\tr(\theta_{E, \chi})$ ne sont rien d'autre que les sommes exponentielles 
\begin{equation*}\label{traces-conj}
S(M, E, \chi)=\sum_{x \in E^\times} \chi(x) \tr(\Frob_{E, x} \hspace{.3mm} | \hspace{.3mm} M)
\end{equation*} que l'on souhaite étudier. On y pense comme aux transformées de Fourier, relatives aux caractères $\chi$, de la \textit{fonction trace} $x \mapsto \tr(\Frob_{E, x} \hspace{.3mm} | \hspace{.3mm} M)$ sur le groupe multiplicatif $E^\times$. Comment varient-elles lorsque le degré des extensions $E\slash k$ tend vers l'infini? 

Voici le résultat principal de cet exposé:

\begin{theo}[Katz]\label{theo:intro-main} Soit $M$ un faisceau pervers dans $\bP(\GG_{m, k})$. Supposons que~$M$ est $\iota$-pur de poids zéro, semisimple, et que les groupes $\mongeo$ et $\monarith$ coïncident. Lorsque le degré des extensions $E \slash k$ tend vers l'infini, les éléments~$\{\theta_{E, \chi}\}_{\chi \hspace{.5mm} \text{bon}}$ s'équirépartissent dans l'espace des classes de conjugaison dans $K$ selon la mesure induite par la mesure de Haar normalisée. 
\end{theo}

En particulier, les sommes $S(M, E, \chi)$ varient comme les traces de matrices aléatoires dans $K$. Il ne reste plus qu'à calculer le groupe $\monarith$. L'avantage du point de vue tannakien est que l'on dispose de maintes techniques pour le faire. On cherche d'abord des bornes par le haut en étudiant la géométrie de $M$ et de ses constructions tensorielles, puis des bornes par le bas en exhibant des éléments explicites dans $\monarith$. Une fois ces contraintes obtenues, la théorie des représentations permet souvent de conclure qu'il n'y a qu'un seul groupe les satisfaisant. Dans le cas des sommes de Gauss, $\mongeo$ est par définition un sous-groupe de $\GL(1)$; comme aucune puissance de convolution de~$M$ n'est égale à l'objet neutre, il s'agit forcément de $\GL(1)$ tout entier, si bien que l'on retrouve le théorème \ref{theo-gauss}. Quant aux objets associés aux sommes d'Evans et de Rudnick, on verra que leurs groupes arithmétiques et géométriques sont dans les deux cas égaux à $\SL(2)$; l'équirépartition  des angles $\theta_{q, \chi}$ selon la mesure de Sato--Tate en découle. 

Le texte est organisé comme suit. Le numéro \ref{sec:rappels-pervers} contient des rappels des résultats sur les faisceaux $\ell$-adiques qui seront utilisés au long du texte; j'ai essayé de les rendre aussi accessibles que possible. Dans le numéro \ref{sec:deligne}, on présente le théorème d'équirépartition  de Deligne et ses applications aux sommes de Kloosterman, puis on explique pourquoi ce résultat permet d'établir la variante additive de l'équirépartition des sommes~$S(M, E, \chi)$. Enfin, le numéro \ref{sec:multplicatif} est consacré aux résultats de la monographie~\cite{katz}, notamment le théorème \ref{theo:intro-main}, une version horizontale et les trois exemples ci-dessus. J'aurais voulu aussi parler des applications du théorème de Katz en théorie analytique des nombres (\cite{BFKMMS},~\cite{HKR}, \cite{katzDirichletI}, \cite{katz-DirichletII}, \cite{KR}, \cite{Xi}), du cas des courbes elliptiques~(\cite{CDRS},~\cite{katz-elliptic}), de l'approche tannakienne aux conjectures d'équirépartition horizontale~\cite{thesis-sawin} et des théorèmes d'annulation générique pour les faisceaux pervers sur les variétés abéliennes (\cite{kramer-weissauer}, \cite{weissauer}), mais cela aurait fait de ce texte un séminaire-fleuve... 

\textit{Remerciements} --- Ils vont tout d'abord à Nick Katz pour la joie des livres orange (renouvelée à chaque lecture!), puis à Henryk Iwaniec et à Emmanuel Kowalski, qui m'ont demandé en premier de présenter ces résultats dans un séminaire informel à l'ETH. Je tiens à remercier aussi Olivier Benoist, Brian Conrad, Nick Katz, Arthur Forey, Emmanuel Kowalski, Corentin Perret\nobreakdash-Gentil, Lior Rosenzweig et Zeév Rudnick pour l'aide qu'ils m'ont apportée au cours de la rédaction de ce texte. 

\section{Rappels sur les faisceaux $\ell$-adiques}\label{sec:rappels-pervers}

Dans ce numéro, on fait un tour rapide de la théorie des systèmes locaux et des faisceaux constructibles $\ell$-adiques. On rappelle notamment la notion de poids et le théorème principal de Weil II, qui sera l'outil essentiel pour obtenir les majorations des sommes exponentielles requises par le critère d'équirépartition  de Weyl. 

\subsection{Groupe fondamental étale}\label{sec:groupe-etale}

Pour tout schéma connexe $X$, notons $\text{\bf{Fét}}_X$ la catégorie des revêtements finis étales~$Y \to X$. Soit~$\xi \colon \Spec(\Omega) \to X$ un \textit{point géométrique} de $X$, où $\Omega$ est un corps algébriquement clos. En associant à chaque $Y \to X$ l'ensemble sous-jacent à la fibre~$Y_\xi=Y \times_X \Spec(\Omega)$, on obtient un foncteur
\begin{displaymath}
\Fib_\xi \colon \text{\bf{Fét}}_X \longrightarrow \text{\bf{Ensembles finis}}.  
\end{displaymath} 

Dans \cite{SGA1}, Grothendieck définit \textit{le groupe fondamental étale avec point base $\xi$} comme le groupe $\pi_1(X, \xi)$ des automorphismes de $\Fib_\xi$, c'est-à-dire les familles $(\sigma_{Y_\xi})_{Y \to X}$ de permutations $\sigma_{Y_\xi}$ des ensembles finis $Y_\xi$ telles que $\sigma_{Y'_\xi}=f \circ \sigma_{Y_\xi}$ pour tout morphisme de~$X$\nobreakdash-schémas $f \colon Y\to Y'$. Par suite, $\pi_1(X, \xi)$ est un groupe \textit{profini}, muni en particulier d'une topologie. Comme en théorie de Galois classique, on~a
\begin{displaymath}
\pi_1(X, \xi) \cong \varprojlim \mathrm{Aut}_X(Y), 
\end{displaymath} la limite étant prise sur les revêtements finis étales \textit{galoisiens}, pour lesquels le schéma~$Y$ est connexe et le groupe $\mathrm{Aut}_X(Y)$ agit transitivement sur $Y_\xi$. 

Puisque $X$ est supposé connexe, les foncteurs fibre $\Fib_{\xi_1}$ et $\Fib_{\xi_2}$ définis par deux points géométriques $\xi_1$ et $\xi_2$ sont isomorphes. Le choix d'un isomorphisme, que l'on appelle parfois \textit{chemin} par analogie avec la topologie, induit un isomorphisme continu $$\pi_1(X, \xi_1) \simeq \pi_1(X, \xi_2)$$ qui ne dépend du chemin qu'à un automorphisme intérieur près. 

De plus, la formation du groupe fondamental est \textit{covariante}: un morphisme de schémas \hbox{$f \colon X \to X'$} induit un homomorphisme continu $f_\ast \colon \pi_1(X, \xi) \to \pi_1(X', f \circ \xi)$. Pour un autre choix de point base $\xi'$ sur~$X'$, on peut encore parler de l'homomorphisme
\begin{displaymath}
f_\ast \colon \pi_1(X, \xi) \to \pi_1(X', \xi'), 
\end{displaymath} mais celui-ci n'a un sens qu'à conjugaison près. 

\begin{exem}\label{fundamental=Galois} Si $X$ est le spectre d'un corps $k$, un point géométrique $\xi$ est un corps algébriquement clos $L$ contenant $k$, et le groupe fondamental $\pi_1(X, \xi)$ s'identifie au groupe de Galois $\Gal(k^{\textup{sep}} \slash k)$ de la clôture séparable $k^{\textup{sep}}$ de $k$ dans~$L$. En particulier, si~$k$ est un corps fini à $q$ éléments, ce groupe est canoniquement isomorphe à la complétion profinie $\widehat{\ZZ}$ de $\ZZ$, avec générateur la substitution de Frobenius $x \mapsto x^q$. On appellera \textit{Frobenius géométrique} l'inverse $F_k$ de ce générateur. Pour une extension~$E$ de degré $n$ sur $k$, l'homomorphisme induit par $\Spec(E) \to \Spec(k)$ envoie $F_E$ sur $F_k^n$. 
\end{exem}

Soient $k$ un corps, $\bar{k}$ une clôture algébrique de $k$, $X$ une variété géométriquement connexe sur $k$ et $X_{\bar{k}}$ son extension des scalaires à $\bar{k}$. Pour tout~$\xi \colon \Spec(\bar{k}) \to X_{\bar{k}}$, on a d'après un théorème de Grothendieck \cite{SGA1} la suite~exacte 
\begin{equation}\label{eqn:suite-exacte-homot}
1 \longrightarrow \pi_1(X_{\bar{k}}, \xi) \longrightarrow \pi_1(X, \xi) \longrightarrow \Gal(\bar{k} \slash k) \longrightarrow 1, 
\end{equation} où les flèches non triviales sont induites par $X_{\bar{k}} \to X$ et $X \to \Spec(k)$. En particulier, si $k$ est un corps fini, on dispose d'une application \textit{degré} $\pi_1(X, \xi) \to \widehat{\ZZ}$. 

Voici la description galoisienne du groupe fondamental d'une courbe. Soit $C$ une courbe projective, lisse et géométriquement connexe sur un corps parfait $k$, avec corps de fonctions $K=k(C)$ et point générique $\eta \colon \Spec(K) \to C$. Fixons une clôture algébrique~$\overline{K}$ de $K$ et notons $\overline{\eta}$ le point géométrique au-dessus de $\eta$ qu'elle définit. Chaque point fermé $x$ de $C$ donne lieu à une valuation discrète $v \colon K^\times \to \ZZ$ de corps résiduel $\kappa(x)$, qui se prolonge en une valuation discrète $\bar{v}$ sur $K^{\textup{sep}}$ (deux tels prolongements sont conjugués sous l'action de $\Gal(K^{\textup{sep}} \slash K)$, et choisir $\bar{v}$ revient à choisir un point géométrique~$\bar{x}$ au-dessus de $x$). Soient $K_{\bar v}^{\textup{sep}}$ et $K_v$ les complétés par rapport aux valuations $\bar v$ et $v$. Le \textit{groupe de décomposition} de $\bar{v}$ est le stabilisateur 
\begin{displaymath}
D_{\bar x}=\{g \in \Gal(K^{\textup{sep}} \slash K) \hspace{.5mm} | \hspace{.5mm} \bar{v} \circ g=\bar{v} \}, 
\end{displaymath} que l'on peut identifier au groupe de Galois de l'extension locale $K_{\bar v}^{\textup{sep}} \slash K_v$. Par passage aux corps résiduels, on obtient une surjection~\hbox{$D_{\bar x} \to \Gal(\bar{k} \slash \kappa(x))$,} où $\bar{k}$ est la clôture algébrique de $k$ dans $\overline{K}$. Son noyau est par définition le \textit{groupe d'inertie} $I_{\bar{x}}$.

\begin{prop}\label{prop:descriptiongalois}
Soit $U$ un ouvert non vide de $C$. Le groupe fondamental étale de~$U$ avec point base $\bar \eta$ est le quotient de $\Gal(K^{\textup{sep}} \slash K)$ par le plus petit sous-groupe fermé distingué contenant $I_{\bar{x}}$ pour tout $x$ dans $U$. 
\end{prop}

En particulier, le groupe fondamental d'un ouvert de la droite projective $\PP^1_k$ s'identifie à un quotient de $\Gal(k(T)^{\sep} \slash k(T));$ c'est le point de vue que l'on a adopté dans l'introduction pour définir les représentations $\cL_\psi$ et $\cL_\chi$.

\subsection{Systèmes locaux $\ell$-adiques et traces de Frobenius}

Dorénavant, $k$ désigne un corps fini à $q$ éléments, $\ell$ un nombre premier différent de la caractéristique de $k$ et $X$ une variété lisse et géométriquement connexe sur $k$, de point générique $\eta$. Pour tout point fermé $x \colon \Spec(E) \to X$, l'homomorphisme continu 
\begin{displaymath}
x_\ast \colon \Gal(\bar{k} \slash E) \longrightarrow \pi_1(X, \bar \eta) 
\end{displaymath} est bien défini à conjugaison près. L'image du Frobenius géométrique $F_E$ par $x_\ast$ est une \textit{classe de conjugaison}\footnote{Avec les notations du paragraphe précédent, quand $X$ est une courbe, il s'agit de l'image de $F_{E}$ par
$\Gal(\bar k \slash E) \to \Gal(\bar k \slash \kappa(x)) \cong D_{\bar x} \slash I_{\bar x} \hookrightarrow \Gal(K^{\textup{sep}} \slash K) \slash I_{\bar{x}} \to \pi_1(X, \bar{\eta})$. Comme l'on a dû choisir un point géométrique $\bar x$ au-dessus de $x$ pour définir $D_{\bar x}$, on obtient seulement une classe de conjugaison.} dans $\pi_1(X, \bar \eta)$, que l'on notera $\Frob_{E, x}$. 

\begin{defi}\label{defi:systemelocal} Un \emph{système local $\ell$-adique} $\cL$ de rang $r$ sur $X$ est la donnée d'un espace vectoriel $\cL_{\bar{\eta}}$ de dimension $r$ sur $\alg_\ell$ et d'un homomorphisme continu 
$$
\rho \colon \pi_1(X, \bar{\eta}) \to \GL(\cL_{\bar\eta}).
$$ 
\end{defi}

Soit $\cL$ un système local $\ell$-adique de rang $r$ sur $X$. Puisque le groupe profini $\pi_1(X, \bar \eta)$ est compact, moyennant le choix d'une base de $\cL_{\bar \eta}$, l'image de $\rho$ est un sous-groupe compact de $\GL(r, \alg_\ell)$. Un tel sous-groupe est toujours contenu dans $\GL(r, E_\lambda)$ pour une extension finie $E_\lambda$ de~$\QQ_\ell$, et l'on peut même choisir la base de sorte que $\rho$ prenne des valeurs dans l'anneau des entiers de cette extension \cite[\S9.0.7]{katz-sarnak}.

L'image de $\Frob_{E, x}$ par la représentation $\rho$ est une classe de conjugaison dans $\GL(\cL_{\bar \eta})$. Sa trace est donc un élément bien défini de $\alg_\ell$, d'où une \textit{fonction trace}
\begin{displaymath}
t_{\cL} \colon X(E) \longrightarrow \alg_\ell, \quad 
x \longmapsto \tr(\Frob_{E, x} \hspace{.3mm} | \hspace{.3mm} \cL)
\end{displaymath} pour toute extension finie $E$ de $k$. Les opérations usuelles sur les représentations (produit tensoriel, dual, etc.) fournissent de nouveaux systèmes locaux dont les fonctions trace s'expriment en termes des anciennes. Par exemple, on a $t_{\cL_1 \otimes \cL_2}=t_{\cL_1}\cdot t_{\cL_2}$ pour le produit tensoriel et, si $f \colon Y \to X$ est un morphisme de $k$-variétés lisses, le système local $f^\ast \cL$ sur $Y$ défini par la représentation $\rho \circ f_\ast$ a fonction trace $t_{f^\ast \cL}=t_{\cL} \circ f$.

\begin{exem} Soit $\alpha$ une unité dans l'anneau des entiers de $\alg_\ell$. L'unique homomorphisme continu $\widehat{\ZZ} \to \alg_\ell^\times$ qui envoie $1$ sur $\alpha$ définit un système local $\ell$-adique de rang un sur $\Spec(k)$; d'après la remarque après la définition \ref{defi:systemelocal}, ils sont tous de cette forme. Par composition avec l'application degré $\pi_1(X, \xi) \to \widehat{\ZZ},$ on obtient un système local~$\ell$\nobreakdash-adique de rang un sur~$X$. Comme l'image de $\Frob_{E, x}$ par la représentation associée est $\alpha^{[E \colon k]}$, il sera noté~$\alpha^{\deg}$. Au vu de la suite exacte~\eqref{eqn:suite-exacte-homot}, les $\alpha^{\deg}$ sont précisément les caractères du groupe fondamental dont la restriction à $\pi_1(X_{\bar k}, \xi)$ est triviale. Si $\alpha=q^n$ pour un entier~$n$, on écrit $\alg_\ell(-n)$ plutôt que $(q^n)^{\deg}$. Si $\alpha \in \alg_\ell$ est une racine carrée de $1/q$, les traces du système local~$\cL \otimes \alpha^{\deg}$, que l'on notera $\cL(1/2)$, sont celles de $\cL$ multipliées par~$1/\sqrt{|E|}$; c'est la \textit{demi-torsion à la Tate} de l'introduction. 
\end{exem}

\begin{exem}\label{exem:lang} Soit $G$ un schéma en groupes commutatifs lisse et connexe sur $k$, la loi de groupe étant notée additivement. Soit $F_G$ le morphisme de Frobenius absolu, c'est\nobreakdash-à\nobreakdash-dire le morphisme de $k$-schémas $F_G \colon G \to G$ qui est l'identité sur l'espace topologique sous-jacent et $x \mapsto x^q$ sur le faisceau structural $\cO_G$. L'\textit{isogénie de Lang} 
\begin{displaymath}
\mathrm{id}-F_G \colon G \to G 
\end{displaymath} est un revêtement fini étale galoisien de groupe $G(k).$ Si $G$ est le groupe additif, c'est le \textit{revêtement d'Artin--Schreier} $x \mapsto x-x^q$ de la droite affine $\IA^1_k$ et, si $G$ est le groupe multiplicatif, c'est le \textit{revêtement de Kummer} $x \mapsto x^{1-q}$ de $\GG_{m, k}$.

Étant donné un caractère $\varphi \colon G(k) \to \alg_\ell^\times$,  la représentation 
\begin{displaymath}
\pi_1(G, \bar \eta) \longrightarrow G(k) \stackrel{\varphi}{\longrightarrow} \alg_\ell^\times
\end{displaymath} obtenue en composant la surjection canonique associée à l'isogénie de Lang avec le caractère $\varphi$ définit un système local $\ell$-adique $\cL_\varphi$ de rang un sur $G$. Plus généralement, pour tout morphisme de $k$-schémas $f \colon X \to G$, on obtient un système local de rang un~$f^\ast \cL_\varphi$ sur $X$ qu'il est coutumier de noter $\cL_{\varphi(f)}$. Pour chaque point de $X$ à valeurs dans~$E$, en raisonnant comme dans le paragraphe \ref{sec:intromethode} de l'introduction, on trouve
\begin{displaymath}
t_{\cL_{\varphi(f)}}(x)=\tr(\Frob_{E, x}\hspace{.5mm} | \hspace{.5mm} \cL_{\varphi(f)})=\varphi(\tr^G_{E \slash k} f(x)),  
\end{displaymath} où $\tr^G_{E \slash k}$ désigne la trace au sens du groupe abélien $G(E)$, c'est-à-dire la trace~$\tr_{E \slash k}$ si~$G$ est le groupe additif, la norme $N_{E \slash k}$ si $G$ est le groupe multiplicatif, etc. 
\end{exem}

\begin{exem} Soit $f \colon \GG_m \to \IA^1$ la fonction $x \mapsto x-x^{-1}$. Le système local $\cL_{\psi_q(f)}$ induit par le caractère additif $\psi_q \colon \FF_q \to \alg^\times$ a trace de Frobenius $\psi_q(x-x^{-1})$ en le point $x \in \FF^\times_q.$ Les sommes d'Evans sont donc égales à 
\begin{displaymath}
S(\chi)=-\sum_{x \in \FF_q^\times} \tr(\Frob_{\FF_q, x} \hspace{.3mm} | \hspace{.3mm} \cL_\chi \otimes \cL_{\psi_q}(1/2)). 
\end{displaymath}
De même, en prenant pour $f \colon \GG_m \setminus \{1\} \to \IA^1$ la fonction $(x+1)\slash (x-1)$, les sommes de Rudnick peuvent se récrire en termes du système local $\cL_{\psi_q(f)}$. 
\end{exem}

Le cadre des systèmes locaux $\ell$-adiques, par opposition aux faisceaux constructibles que l'on introduira dans le paragraphe suivant, est peu adapté aux méthodes cohomologiques. Par exemple, si $\chi \colon k^\times \to \alg_\ell^\times$ est un caractère et $f \colon X \to \GG_{m, k}$ un morphisme de $k$-variétés, $\chi(f(x))$ est bien la trace de Frobenius en $x$ du système local $\cL_{\chi(f)}$ sur~$X$, mais il est souvent plus commode d'utiliser l'égalité 
\begin{displaymath}
\sum_{x \in X(k)} \chi(f(x))=\sum_{t \in k^\times} \chi(t) |f^{-1}(t)|
\end{displaymath} et de voir $t \mapsto |f^{-1}(t)|$ comme une fonction trace sur $\GG_{m, k}$. Ceci correspond à étudier la variation en $t$ de la cohomologie des fibres $f^{-1}(t)$ et l'on doit s'attendre à ce que seule la restriction à un ouvert de lissité de $f$ soit la trace d'un système local. 

\subsection{Faisceaux constructibles, poids et cohomologie}\label{faisceaux-constr}

Vers la fin des années 50, Grothendieck définit une \og topologie \fg{} sur les schémas où la notion classique d'ouvert d'une variété $X$ est remplacée par les morphismes étales vers $X$. Si $A$ est un anneau de torsion, par exemple $\ZZ/\ell^n\ZZ$, on peut parler de faisceaux localement constants (aussi appelés \textit{lisses}) et de faisceaux constructibles de $A$-modules sur $X$, ces derniers ayant la propriété qu'il existe une stratification de $X$ telle que la restriction à chaque strate soit un faisceau lisse. Par un double passage à la limite, d'abord sur les entiers $n$, puis sur les extensions finies de $\QQ_\ell$, on dégage la notion de~$\alg_\ell$\nobreakdash-faisceau constructible. On ne se contentera ici que de quelques aperçus.

Soit $\cF$ un faisceau constructible sur $X$. La fibre $\cF_\xi$ en un point géométrique $\xi$ est un espace vectoriel de dimension finie sur $\alg_\ell$. Si $\cF$ est lisse, le groupe fondamental~$\pi_1(X, \bar \eta)$ agit continument sur la fibre~$\cF_{\bar \eta}$ en un point géométrique générique $\bar \eta$ et le foncteur~\hbox{$\cF \mapsto \cF_{\bar \eta}$} induit une équivalence de catégories entre les~$\alg_\ell$\nobreakdash-faisceaux lisses et les systèmes locaux $\ell$-adiques sur $X$. En général, quel que soit~$\cF$, la fibre $\cF_{\bar x}$ en un point géométrique $\bar{x}$ au-dessus de~$x \in X(E)$ est un système local $\ell$-adique sur $\Spec(E)$, d'où une action du Frobenius $F_E$ qui permet d'étendre la définition de fonction trace à tout faisceau constructible. 

Soit $C$ une courbe projective, lisse et géométriquement connexe sur $k$. Si $\cF$ est un faisceau lisse sur un ouvert non vide $V \subseteq C$, pour tout point fermé $x \in C$, le groupe de décomposition $D_{\bar x}$ agit sur $\cF_{\bar \eta}$ au travers de la flèche $D_{\bar x} \to \pi_1(V, \bar \eta)$. On en déduit une action du groupe d'inertie $I_{\bar x}$ pour tout $x \in C \setminus V$, celle de $I_{\bar x}$ pour $x \in V$ étant triviale par la description galoisienne du groupe fondamental (proposition \ref{prop:descriptiongalois}). Concrètement, un $\alg_\ell$-faisceau constructible $\cF$ sur une courbe $U$ est la donnée:  
\begin{enumerate}
\item[$\bullet$] d'un ouvert non vide $V \hookrightarrow U$ et d'un faisceau $\ell$-adique lisse $\cF_{\bar \eta}$ sur $V$ correspondant à une représentation continue $\pi_1(V, \bar{\eta}) \to \GL(\cF_{\bar{\eta}}),$
\item[$\bullet$] pour chaque point fermé $u \in U \setminus V$, d'une représentation continue de $\Gal(\bar{k} \slash \kappa(u))$ dans un $\alg_\ell$-espace vectoriel de dimension finie $\cF_{\bar{u}},$ 
\item[$\bullet$] des \textit{flèches de spécialisation} $\Gal(\bar{k} \slash \kappa(u))$-équivariantes 
$\mathrm{sp}_{\bar u, \bar \eta} \colon \cF_{\bar u} \to \cF_{\bar \eta}^{I_{\bar u}}$ (les données de recollement).  
\end{enumerate} 

\begin{exem}\label{exemp:imagedirecte}
Quelques opérations naturelles sur les faisceaux constructibles, par exemple les images directes par une immersion ouverte, admettent aussi une description en ces termes. Si $\cF$ est un faisceau constructible sur $U$ et $j \colon U \hookrightarrow U'$ une immersion ouverte, le \textit{prolongement par zéro} $j_! \cF$ (resp. l'\textit{image directe} $j_\ast \cF$) est le faisceau constructible sur~$U'$ qui coïncide avec $\cF$ sur~$U$ et qui a fibre $(j_! \cF)_{\bar x}=0$ (resp.~$(j_\ast \cF)_{\bar{x}}=\cF^{I_{\bar x}}_{\bar \eta}$, avec flèche de spécialisation l'identité) en tout point fermé $x$ dans le complémentaire. L'\textit{image directe supérieure} $R^1j_\ast \cF$ est le faisceau constructible supporté dans $U' \setminus U$ et ayant pour fibre en un point $x$ les coinvariants sous l'inertie $(R^1j_\ast \cF)_x=(\cF_{\bar \eta})_{I_{\bar x}}(-1)$.  
\end{exem}

Soit $X$ une variété lisse et géométriquement connexe, purement de dimension $d$, sur un corps fini~$k$. À un $\alg_\ell$-faisceau constructible $\cF$ sur $X$ on associe des \textit{groupes de cohomologie} $H^i(X_{\bar k}, \cF)$ et de \textit{cohomologie à support compact} $H^i_c(X_{\bar k}, \cF)$. Ce sont des espaces vectoriels de dimension finie sur $\alg_\ell$, munis d'une action continue de $\Gal(\bar k \slash k)$, qui s'annulent en degrés $i<0$ et $i>2d$. Pour chaque $i$, il y a une flèche canonique d'\og oubli des supports \fg{} $H^i_c(X_{\bar k}, \cF) \to H^i(X_{\bar k}, \cF)$ qui est un isomorphisme lorsque $X$ est propre. Si le faisceau $\cF$ est lisse sur $X$, on a
\begin{displaymath}
H^0(X_{\bar k}, \cF)=(\cF_{\bar \eta})^{\pi_1(X_{\bar k}, \bar \eta)}, \qquad H^{2d}_c(X_{\bar k}, \cF)=(\cF_{\bar \eta})_{\pi_1(X_{\bar k}, \bar \eta)}(-d)
\end{displaymath} et le groupe de Galois opère sur ces espaces par le biais de l'action de $\pi_1(X, \bar \eta)$ sur $\cF_{\bar \eta}$ et de l'identification $\Gal(\bar k \slash k)=\pi(X, \bar \eta) \slash \pi_1(X_{\bar k}, \bar \eta)$ donnée par la suite exacte \eqref{eqn:suite-exacte-homot}. En général, un faisceau constructible~$\cF$ peut avoir des sections globales ponctuelles non nulles (elles correspondent aux éléments dans le noyau des flèches de spécialisation).  

 Le lien entre les sommes exponentielles et l'action du Frobenius géométrique $F_k$ sur la cohomologie à support compact est fourni par la formule des traces \cite{lefschetz}: 

\begin{theo}[Formules des traces de Grothendieck] Pour toute extension finie~$E$ de $k$, on a l'égalité
\begin{displaymath}
\sum_{x \in X(E)} \tr(\Frob_{E, x} \hspace{.5mm} | \hspace{.5mm} \cF)=\sum_{i=0}^{2d} (-1)^i \tr(F_E \hspace{.5mm} | \hspace{.5mm} H^i_c(X_{\bar{k}}, \cF)). 
\end{displaymath} 
\end{theo}

Cette formule n'acquiert toute sa force qu'en combinaison avec des estimées sur les valeurs propres de l'endomorphisme $F_E$: le formalisme des poids de la cohomologie. Pour pouvoir parler de la taille de ces valeurs propres, qui sont des éléments dans $\alg_\ell$, on fixe un plongement $\iota \colon \alg_\ell \hookrightarrow \CC$. 

\begin{defi}\label{defi-poids} Soient $w$ un entier et $\cF$ un faisceau constructible sur $X$. 
\begin{enumerate}[wide =1em]
\item[\rm (a)] On dit que $\cF$ est {\rm{$\iota$-pur de poids $w$}} si, pour toute extension finie $E$ de $k$ et pour tout point $x \in X(E)$, les valeurs propres de $\Frob_{E, x}$ agissant sur la fibre $\cF_{\bar x}$, où $\bar{x}$ est un point géométrique quelconque au-dessus de $x$, ont valeur absolue~$|E|^{w/2}$ quand on les regarde dans $\CC$ à travers le plongement $\iota$. 

\item[\rm (b)] On dit que $\cF$ est {\rm $\iota$-mixte} s'il admet une filtration finie croissante par des sous\nobreakdash-faisceaux constructibles $(\cF_i)_{i \in \ZZ}$ dont les gradués $\gr_i=\cF_i \slash \cF_{i-1}$ sont $\iota$-purs de poids $i$. Si les poids des $\gr_i$ sont majorés (resp. minorés) par $w$, on dira que $\cF$ est~$\iota$\nobreakdash-mixte de poids $\leq w$ (resp. $\geq w$). 
\end{enumerate}
\end{defi}

Le poids d'un faisceau dépend a priori fortement du plongement $\iota$; on dira que  $\cF$ est \textit{pur de poids $w$} s'il est $\iota$-pur de poids $w$ pour tout $\iota \colon \alg_\ell 
\hookrightarrow \CC$. Pour des questions concernant les faisceaux purs, quitte à les demi-tordre à la Tate, on pourra toujours supposer que le poids est zéro. Le théorème principal de Weil II relie le poids d'un faisceau constructible au sens de la définition ci-dessus aux valeurs propres de Frobenius agissant sur sa cohomologie et sa cohomologie à support compact.

\begin{theo}[Deligne, \cite{weil2}]\label{theo-poids} Soit $\cF$ un faisceau constructible sur $X$. Si $\cF$ est~$\iota$\nobreakdash-mixte de poids $\leq w$, alors chaque $H^i(X_{\bar k}, \cF)$ est $\iota$-mixte de poids $\geq w+i$ et chaque $H^i_c(X_{\bar k}, \cF)$ est $\iota$-mixte de poids $\leq w+i$. En particulier, si la flèche d'oubli des supports est un isomorphisme, $H^i(X_{\bar k}, \cF) \cong H^i_c(X_{\bar k}, \cF)$ est $\iota$-pur de poids $w+i$. 
\end{theo}

On renvoie aux travaux de Laumon \cite{laumon} et de Katz \cite{katz4lectures} pour des simplifications de la preuve originale, basées sur l'incarnation faisceautique de la transformation de Fourier. Voici un exemple frappant d'application de ces résultats qui reviendra dans la preuve du théorème d'équirépartition  de Deligne: 

\begin{exem}[Inégalités de Lang--Weil, \cite{LangWeil}]\label{exmp:inegalitesLangWeil} En prenant pour $\cF$ le faisceau constant~$\alg_\ell$, on trouve que le nombre de points $E$-rationnels de $X$ est donné par 
\begin{equation*}\label{eqn:numberofpoints}
|X(E)|=\sum_{i=0}^{2d} (-1)^i \tr(F_E \hspace{.3mm} | \hspace{.3mm} H^i_c(X_{\bar k}, \alg_\ell)). 
\end{equation*} 

Comme $\alg_\ell$ est pur de poids zéro, $H^i_c(X_{\bar k}, \alg_\ell)$ est mixte de poids $\leq i$ d'après le théorème \ref{theo-poids}, ce qui signifie que les valeurs propres de $F_E$ agissant sur cet espace ont toutes valeur absolue au plus $|E|^{i/2}$ quel que soit le plongement de $\alg_\ell$ dans $\CC$. Posons $$
b^i_c(X)=\dim_{\alg_\ell} H^i_c(X_{\bar k}, \alg_\ell), \qquad A(X)=\sum_{i=0}^{2d-1} b^i_c(X).
$$ On a vu que $H^{2d}_c(X_{\bar k}, \alg_\ell)$ est le $\alg_\ell$-espace vectoriel de dimension un sur lequel $F_E$ agit par multiplication par $|E|^d$, d'où les inégalités 
\begin{equation}\label{eqn:LangWeil}
\left| |X(E)|-|E|^d \right| \leq \sum_{i=0}^{2d-1} b^i_c(X) |E|^{i/2} \leq A(X) |E|^{(2d-1)/2}. 
\end{equation} En particulier, la variété $X$ a un point $E$-rationnel dès que $|E|>A(X)^2$. 
\end{exem}

\subsection{Un peu de ramification}\label{sec-ramification}

Dans ce qui suit, on rappelle une formule due à Grothendieck--Ogg--Shafarevich qui exprime la caractéristique d'Euler d'un faisceau constructible sur une courbe en termes de sa ramification à l'infini. Soit $C$ une courbe projective, lisse et géométriquement connexe sur un corps parfait~$k$ de caractéristique $p>0$. Pour chaque point fermé $x$ de~$C$ et chaque point géométrique~$\bar x$ au-dessus de $x$, nous avons une suite exacte  
\begin{displaymath}
1 \longrightarrow I_{\bar x} \longrightarrow D_{\bar x} \longrightarrow \Gal(\bar{k} \slash k) \longrightarrow 1.  
\end{displaymath} 

Le groupe d'inertie $I_{\bar x}$ contient un unique $p$-sous-groupe de Sylow $P_{\bar x}$, sa partie \textit{sauvage}, et le quotient $I_{\bar x}^{\mathrm{mod}}=I_{\bar x} \slash P_{\bar x}$ (l'inertie \textit{moderée}) est canoniquement isomorphe~à
\begin{displaymath}
\widehat{\ZZ}'(1)=\varprojlim_{(N, p)=1} \mu_N({\bar k})=\varprojlim_{E/k \text{ finie}} E^\times.
\end{displaymath} Il s'agit du premier cran d'une filtration décroissante $(I_{\bar x}^{(r)})_{r \in \QQ_{\geq 0}}$ par des sous-groupes fermés distingués de $I_{\bar x}$ que l'on appelle la \textit{filtration de ramification en numérotation supérieure} \cite[Ch. IV]{corps-locaux}: on a $I_{\bar x}^{(0)}=I_{\bar x}$, la partie sauvage $P_{\bar x}$ est l'adhérence de $\bigcup_{r>0} I_{\bar x}^{(r)}$ dans $I_{\bar x}$, et $I^{(r)}_{\bar x}=\bigcap_{0<s<r} I^{(s)}_{\bar x}$ pour tout $r \in \QQ_{>0}$. 

Soit $I_{\bar x} \to \GL(W)$ une représentation continue dans un espace vectoriel $W$ de dimension finie sur $\alg_\ell$. On dit que $W$ est \textit{modérée} si $P_{\bar x}$ agit trivialement et \textit{totalement sauvage} si $W^{P_{\bar x}}=0$. La représentation $W$ admet une décomposition unique
\begin{equation}\label{eqn:dec-swan}
W=\bigoplus_{r \geq 0} W(r),
\end{equation} indexée par des nombres rationnels, en des sous-espaces vectoriels $W(r) \subseteq W$ stables sous l'action de $P_{\bar x}$ et vérifiant les propriétés suivantes: 
\begin{displaymath}
W(0)=W^{P_{\bar x}}, \qquad W(r)^{I_{\bar x}^{(r)}}=0 \ \text{ si }r>0, \qquad W(r)^{I_{\bar x}^{(r')}}=W(r) \ \text{ si } r'>r.  
\end{displaymath} Cette décomposition ne comporte qu'un nombre fini de valeurs de $r$ pour lesquelles~$W(r)$ n'est pas nul, que l'on appellera les \textit{sauts}. De plus, $r \dim W(r)$ est un entier positif. On définit le \textit{conducteur de Swan} de la représentation $W$ comme l'entier positif
\begin{displaymath}
\mathrm{Sw}(W)=\sum_{r \geq 0} r \dim W(r). 
\end{displaymath}

Soient $V$ un ouvert non vide de $C$ et $\cF$ un faisceau lisse sur $V$, vu comme représentation $\ell$-adique de $\pi_1(V, \bar \eta)$. Par restriction, on en déduit une représentation du groupe d'inertie $I_{\bar x}$ pour tout point fermé $x \in C \setminus V$. On notera $\mathrm{Sw}_{\bar x}(\cF)$ son conducteur de Swan. Il s'agit d'un entier positif qui est nul si et seulement si la représentation de~$I_{\bar x}$ est modérée, auquel cas on dira que $\cF$ a \textit{ramification modérée} en $x$. Si $V=\GG_{m, k}$, les systèmes locaux modérément ramifiés à l'infini sont les représentations du groupe fondamental modéré $\pi_1^{\textup{mod}}(\GG_{m, k})$, qui s'insère dans une suite exacte  
\begin{displaymath}
1 \longrightarrow \widehat{\ZZ}'(1) \longrightarrow \pi_1^{\textup{mod}}(\GG_{m, k}) \longrightarrow \Gal(\bar k \slash k) \longrightarrow 1. 
\end{displaymath}

\begin{exem} Soit $k$ un corps fini. Le faisceau d'Artin--Schreier $\cL_\psi$ sur $\IA^1_k$ a ramification totalement sauvage à l'infini et son conducteur de Swan vaut $1$. Le faisceau de Kummer $\cL_\chi$ sur $\GG_{m, k}$ a ramification modérée en $0$ et $\infty$. Par extension, on appelera \textit{faisceau de Kummer} tout système local~$\ell$\nobreakdash-adique de rang un sur $\mathbb{G}_{m, k}$ correspondant à un caractère de $\pi_1^{\textup{mod}}(\GG_{m, k})$.
\end{exem}

Soient maintenant $U \subseteq C$ un ouvert contenant $V$ et $\cF$ un faisceau constructible sur~$U$. Pour tout point fermé $u \in U \setminus V$, on définit l'entier
$$
\dr_{\bar u}(\cF)=\dim(\cF_{\bar \eta}) - \dim(\cF_{\bar{u}}). 
$$ Si $\cF$ n'a pas de sections globales ponctuelles non nulles, alors toutes les flèches de spécialisation $\cF_{\bar{u}} \to \cF_{\bar \eta}^{I_{\bar u}}$ sont injectives; dans ce cas, $\dr_{\bar u}(\cF)$ est positif et le faisceau~$\cF$ est lisse en $\bar u$ si et seulement si $\dr_{\bar u}(\cF)=0$. 

La formule de Grothendieck--Ogg--Shafarevich, exposée dans ce séminaire par Raynaud \cite{raynaud}, est l'énoncé que la caractéristique d'Euler de $\cF$ s'exprime ainsi:  

\begin{theo}[Grothendieck--Ogg--Shafarevich] Soient $\cF$ un faisceau $\ell$-adique constructible sur $U$ et $V \subseteq U$ un ouvert non vide sur lequel $\cF$ est lisse. Alors
\begin{displaymath}
\chi(U_{\bar{k}}, \cF)=\mathrm{rg}(\cF) \chi(U_{\bar{k}})-\sum_{x \in (C \setminus V)(\bar k)}  \Sw_{\bar x}(\cF)-\sum_{u \in (U \setminus V)(\bar k)} \dr_{\bar u}(\cF). 
\end{displaymath}
\end{theo}

\subsection{Catégorie dérivée et faisceaux pervers}

Dans ce paragraphe, $k$ désigne ou bien un corps fini, ou bien la clôture algébrique d'un corps fini. Soit $\pi \colon X \to \Spec(k)$ une variété sur $k$. Grothendieck a défini la \textit{catégorie dérivée}~$D^b_c(X, \alg_\ell)$ des complexes bornés de $\alg_\ell$-faisceaux sur~$X$ à cohomologie constructible et l'a munie d'un formalisme des six opérations:  
\begin{itemize}[wide =1em]
\item[$\bullet$] Pour tout morphisme $f \colon X \to Y$, on dispose des foncteurs \textit{image directe} et \textit{image directe à support compact}
\begin{displaymath}
Rf_\ast, Rf_! \colon D^b_c(X, \alg_\ell) \longrightarrow D^b_c(Y, \alg_\ell)
\end{displaymath} et des foncteurs \textit{image inverse} et \textit{image inverse exceptionnelle}
\begin{displaymath}
f^\ast, f^! \colon D^b_c(Y, \alg_\ell) \longrightarrow D^b_c(X, \alg_\ell). 
\end{displaymath} En particulier, on obtient pour tout objet~$M$ dans $D^b_c(X, \alg_\ell)$ des complexes de~$\alg_\ell$\nobreakdash-espaces vectoriels $\mathrm{R}\pi_\ast M$ et~$\mathrm{R}\pi_! M$ que l'on notera~$\mathrm{R}\Gamma(X, M)$ et~$\mathrm{R}\Gamma_c(X, M).$ Si~$M$ est de la forme $\cF[n]$, la cohomologie en degré $i$ de ces complexes s'identifie à la cohomologie $H^{i+n}(X_{\bar k}, \cF)$ et à la cohomologie à support compact $H^{i+n}_c(X_{\bar k}, \cF)$. 

\item[$\bullet$] On dispose d'un bifoncteur \textit{produit tensoriel} 
\begin{displaymath}
-\otimes- \colon D^b_c(X, \alg_\ell) \times D^b_c(X, \alg_\ell) \longrightarrow D^b_c(X, \alg_\ell)
\end{displaymath} et d'un bifoncteur \textit{Hom interne} 
\begin{displaymath}
R\mathcal{H}om(-, -) \colon D^b_c(X, \alg_\ell)^{\textup{op}} \times D^b_c(X, \alg_\ell) \longrightarrow D^b_c(X, \alg_\ell). 
\end{displaymath} Le \textit{produit extérieur} des objets $M$ dans $D^b_c(X, \alg_\ell)$ et $N$ dans $D_b(Y, \alg_\ell)$ est l'objet 
\begin{displaymath}
M \boxtimes N=\pr_X^\ast M \otimes \pr_Y^\ast N
\end{displaymath} dans $D^b_c(X \times Y, \alg_\ell)$, où $\pr_X$ et~$\pr_Y$ désignent les projections de $X \times Y$ sur $X$ et $Y$.

\item[$\bullet$] On définit le \textit{complexe dualisant} comme $\omega_X=\pi^! \alg_\ell$ et le foncteur de \textit{dualité de Verdier} $\DD \colon D^b_c(X, \alg_\ell)^{\textup{op}}\to D^b_c(X, \alg_\ell)$ comme
\begin{displaymath}
\DD(M)=R\mathcal{H}om(M, \pi^! \alg_\ell). 
\end{displaymath} La dualité échange $Rf_\ast$ et $Rf_!$ (resp. $f^\ast$ et $f^!$) pour tout morphisme $f \colon X \to Y$. Si $X$ est lisse et purement de dimension $d$, le complexe dualisant est égal à \hbox{$\omega_X=\alg_\ell(d)[2d]$} et l'on retrouve la dualité de Poincaré~$H^i(X, \alg_\ell) \cong H^{2d-i}_c(X, \alg_\ell)^\ast(-d)$ comme cas particulier. Plus généralement,~$f^!=f^\ast(d)[2d]$ si $f$ est lisse de dimension relative $d$. 
\end{itemize}

Dans le cas où $k$ est un corps fini, on définit la fonction trace d'un objet $M$ de la catégorie dérivée comme la somme alternée $t_{M}=\sum (-1)^i t_{\cH^i(M)}$ des fonctions trace de ses faisceaux de cohomologie. La formule des traces de Grothendieck est alors l'égalité 
\begin{equation}\label{eq:formuledestraces}
t_{Rf_!M}(y)=\sum_{f(x)=y} t_{M}(x). 
\end{equation}
De même, la notion de poids s'étend à la catégorie dérivée comme suit: un objet~$M$ dans~$D^b_c(X, \alg_\ell)$ est dit \textit{$\iota$-mixte de poids $\leq w$} si tous ses faisceaux de cohomologie $\cH^i(M)$ sont $\iota$-mixtes de poids $\leq w+i$, au sens de la définition \ref{defi-poids}, et \textit{$\iota$-mixte de poids} $\geq w$ si son dual de Verdier $\DD(M)$ est $\iota$-mixte de poids $\leq -w$. On dit enfin que $M$ est \textit{$\iota$-pur de poids $w$} s'il est en même temps $\iota$-mixte de poids $\leq w$ et $\geq w$. Le théorème principal de Weil II peut alors se formuler ainsi:  

\begin{theo}[Deligne]\label{theo:del-cat} Si $f \colon X \to Y$ est un morphisme de $k$-variétés, le foncteur $Rf_!$ envoie les objets $\iota$-mixtes de poids $\leq w$ sur des objets $\iota$-mixtes de poids $\leq w$. Par dualité de Verdier, le foncteur $Rf_\ast$ envoie les objets $\iota$-mixtes de poids $\geq w$ sur des objets $\iota$-mixtes de poids $\geq w$. 
\end{theo}

En s'inspirant des travaux de Goresky et MacPherson sur l'homologie d'intersection d'un espace topologique singulier, Beilinson, Bernstein, Deligne et Gabber ont introduit la notion de faisceau pervers dans \cite{bbd}.

\begin{defi}
Un objet $M$ de $D^b_c(X, \alg_\ell)$ est dit {\rm{semipervers}} si tous ses faisceaux de cohomologie $\cH^i(M)$ satisfont à la condition de support 
\begin{displaymath}
\dim \mathrm{supp}(\cH^i(M)) \leq -i.  
\end{displaymath} On appelle $M$ un {\rm{faisceau pervers}} si tant $M$ que $\DD(M)$ sont semipervers. 
\end{defi}

Par exemple, si $X$ est lisse de dimension $d$ et $\cL$ est un système local $\ell$-adique sur $X$, l'objet $\cL[d]$ est un faisceau pervers car la condition de support est vérifiée trivialement et son dual $\DD(\cL[d])=\cL^\vee(d)[d]$ est encore de la même forme. 

La sous-catégorie pleine $\Perv(X, \alg_\ell)$ de $D^b_c(X, \alg_\ell)$ formée des faisceaux pervers est abélienne et tout objet est extension successive d'un nombre fini d'objets simples \cite{bbd}. Qui plus est, la catégorie dérivée de complexes bornés dans~$\Perv(X, \alg_\ell)$ est équivalente à $D^b_c(X, \alg_\ell)$, d'où des foncteurs de \textit{cohomologie perverse} 
\begin{displaymath}
\tensor[^p]{\cH}{^0} \colon D^b_c(X, \alg_\ell) \longrightarrow \Perv(X, \alg_\ell)
\end{displaymath} et $\tensor[^p]{\cH}{^i}=\tensor[^p]{\cH}{^0}[i]$ pour tout entier $i$. Si $f \colon X \to Y$ est un morphisme affine, alors $Rf_\ast$ préserve la semiperversité (théorème d'annulation d'Artin) et si~$f \colon X \to Y$ est un morphisme lisse de dimension relative $d$, le foncteur $f^\ast[d]$ préserve la perversité. 

La plupart du temps nous travaillerons avec des faisceaux pervers sur une courbe, auquel cas les objets admettent une description plus élémentaire. 

\begin{exem}[Faisceaux pervers sur les courbes]\label{perv-courbes}
Soit $X$ une courbe lisse et géométriquement connexe sur un corps fini $k$. Un faisceau pervers sur $X$ est un complexe de faisceaux constructibles $M$ dans~$D^b_c(X, \alg_\ell)$ satisfaisant aux conditions suivantes: $\cH^{i}(M)=0$ pour tout~$i \notin \{-1, 0\}$, le faisceau constructible $\cH^{-1}(M)$ n'a pas de sections globales ponctuelles non nulles, et le faisceau constructible $\cH^0(M)$ est supporté dans des points. En posant $M_{\pct}=\cH^0(M)$ et $M_{\npct}=\cH^{-1}(M)[1]$, on trouve une suite exacte  
\begin{equation}\label{eqn:devis}
0 \longrightarrow M_{\npct} \longrightarrow M \longrightarrow M_{\pct} \longrightarrow 0.  
\end{equation}
De plus, les objets simples de la catégorie $\Perv(X, \alg_\ell)$ sont d'un de ces deux types: 
\begin{itemize}[wide =1em]
\item[$\bullet$] \textit{ponctuels}: de la forme $\alpha^{\deg} \otimes \delta_x$ pour un point fermé $x$ de $X$ et une unité $\ell$-adique~$\alpha$; 
\item[$\bullet$] \textit{non ponctuels}: de la forme $(j_\ast \cL)[1]$, où $j \colon U \hookrightarrow X$ est l'inclusion d'un ouvert non vide et $\cL$ est un système local $\ell$-adique simple sur $U$. On les appelle \textit{extensions intermédiaires}. 
\end{itemize}
\end{exem}

\section{Le théorème d'équirépartition  de Deligne}\label{sec:deligne}

Dans ce numéro, on expose le théorème d'équirépartition de Deligne \cite[\S 3.5]{weil2} suivant les variantes qu'en ont données Katz \cite[Ch.\,3]{katz-kloos} et Katz--Sarnak \cite[\S 9.2]{katz-sarnak}. Sous des hypothèses assez faibles, il affirme que le groupe de monodromie d'un système local~$\ell$\nobreakdash-adique gouverne la répartition de ses traces de Frobenius. Ces groupes sont en général difficiles à calculer, mais Katz a réussi à les déterminer complètement pour les faisceaux associés aux sommes de Kloosterman; combiné avec le résultat de Deligne, ceci permet d'obtenir le théorème \ref{thm:equi-kloos} dans l'introduction, ainsi que de l'étendre à un énoncé d'équirépartition pour les sommes $\Kl_n(a, q)$ en plusieurs variables. Dans le dernier paragraphe, on analyse la répartition des traces de la transformée de Fourier d'un faisceau pervers sur la droite affine à l'aide du théorème de Deligne, puis on explique pourquoi cette approche est vouée à l'échec pour le groupe multiplicatif. 

\subsection{Préliminaires sur les représentations}

Dans le théorème d'équirépartition  de Deligne, les groupes de monodromie sont des groupes algébriques réductifs $G$ sur $\alg_\ell$ et l'on définit des classes de conjugaison dans des sous-groupes compacts maximaux $K$ de leurs points complexes. Pour appliquer le critère d'équirépartition  de Weyl, on aura besoin de relier les représentations continues de $K$ à des systèmes locaux $\ell$-adiques, vus comme des $\alg_\ell$-représentations algébriques~de~$G$. 

\begin{theo}\label{theo:equiv} Soit $G$ un groupe réductif (pas nécessairement connexe) sur $\alg_\ell$. Pour chaque choix d'un plongement $\iota \colon \alg_\ell \hookrightarrow \CC$ et d'un sous-groupe compact maximal~$K$ de $G(\CC)$, les foncteurs extension des scalaires, évaluation dans les points complexes et restriction à $K$ induisent des bijections entre les classes d'isomorphisme de
\begin{enumerate}
\item[\rm (a)] $\alg_\ell$-représentations du $\alg_\ell$-groupe algébrique $G$,  
\item[\rm (b)] $\CC$-représentations du $\CC$-groupe algébrique $G_{\CC}$,  
\item[\rm (c)] représentations holomorphes du groupe de Lie complexe~$G(\CC)$, 
\item[\rm (d)] représentations continues du groupe compact $K$, 
\end{enumerate} toutes les représentations étant supposées de dimension finie. 
\end{theo}

\begin{proof} Les catégories (b), (c) et (d) sont en fait équivalentes. Pour démontrer l'équivalence entre (b) et (c), on se ramène facilement au cas où $G$ est connexe, qui est traité dans \cite[Prop.\,D.2.1]{conrad}. L'équivalence entre (c) et (d), sans des hypothèses de connexité, est le contenu de la proposition D.3.2 et de l'exemple~D.3.3 dans \textit{loc.\,cit.} Quant à~(a), on a tautologiquement une équivalence de catégories si le plongement $\iota$ est un isomorphisme de corps. En général, étant donnée une représentation complexe $\rho' \colon G_{\CC} \to \GL(V')$, il existe un espace vectoriel~$V$ sur~$\alg_\ell$ et une représentation~$\rho \colon G \to \GL(V)$, unique à isomorphisme près, telle que~$\rho' \simeq \rho_{\CC}$.  
\end{proof}

\subsection{Le théorème d'équirépartition  de Deligne}\label{sec:theodeligne}

Soient $k$ un corps fini à $q$ éléments, $\bar{k}$ une clôture algébrique de $k$, $X$ une variété lisse et géométriquement connexe sur $k$ et $X_{\bar k}$ la variété sur $\bar k$ qui s'en déduit par extension des scalaires. Fixons un nombre premier $\ell$ distinct de la caractéristique de~$k$, une clôture algébrique~$\alg_\ell$ de $\QQ_\ell$ et un plongement $\iota \colon \alg_\ell \hookrightarrow \CC$. Soit $\cL$ un faisceau~$\ell$\nobreakdash-adique lisse de rang $r$ sur $X$, correspondant à une représentation continue
\begin{displaymath}
\rho \colon \pi_1(X, \bar \eta) \to \GL(r, \alg_\ell).
\end{displaymath} 

\begin{defi} Les groupes de monodromie arithmétique $G_{\arith, \hspace{.3mm} \cL}$ et géométrique $G_{\geom, \hspace{.3mm} \cL}$ de $\cL$ sont les groupes algébriques sur $\alg_\ell$ obtenus comme l'adhérence de Zariski de $\rho(\pi_1(X, \bar \eta))$ et de $\rho(\pi_1(X_{\bar k}, \bar \eta))$ dans $\GL(r, \alg_\ell).$ 
\end{defi}

En particulier, $G_{\geom, \hspace{.3mm} \cL}$ est un sous-groupe de $G_{\arith, \hspace{.3mm} \cL}.$ Comme les classes de conjugaison de Frobenius sont denses dans le groupe fondamental par une forme du théorème de Cebotarev, $G_{\arith, \hspace{.3mm} \cL}$ est aussi le plus petit sous-groupe algébrique contenant leurs images par la représentation. Dans Weil~II, Deligne démontre que, si $\cL$ est $\iota$-pur d'un certain poids, la représentation de $\pi_1(X_{\bar k}, \bar \eta)$ associée est complètement réductible~\hbox{\cite[Thm.\,3.4.1\,(iii)]{weil2}} et le groupe de monodromie géométrique $G_{\geom, \hspace{.3mm} \cL}$ est donc réductif. Qui plus est, la composante neutre du groupe de monodromie géométrique~$G_{\geom, \hspace{.3mm} \cL}$ est semisimple d'après \cite[Cor.\,1.3.9]{weil2}. 

Supposons que $\cL$ est $\iota$-pur de poids zéro et $G_{\geom, \hspace{.3mm} \cL}=G_{\arith, \hspace{.3mm} \cL}$, de sorte que la représentation $\rho$ se factorise par~$G_{\geom, \hspace{.3mm} \cL}(\alg_\ell)$. Choisissons un sous\nobreakdash-groupe compact maximal~$K$ de $G_{\geom, \hspace{.3mm} \cL}(\CC)$ et notons $K^\#$ l'ensemble de ses classes de conjugaison. Le but est d'associer à chaque point de $X$ à valeurs dans une extension finie~$E$ de $k$ un élément de $K^\#$ dont la trace soit égale à~$\tr(\Frob_{E, x} \hspace{.3mm} | \hspace{.3mm} \cL)$. Soient $x$ un tel point et $\Frob_{E, x}$ la classe de conjugaison qu'il définit dans $\pi_1(X, \bar \eta)$. Son image par la représentation $\rho$ est une classe de conjugaison 
\begin{displaymath}
\rho(\Frob_{E, x}) \in G_{\geom, \hspace{.3mm} \cL}(\alg_\ell)^\# \subset G_{\geom, \hspace{.3mm} \cL}(\CC)^\#. 
\end{displaymath} Soit $\rho(\Frob_{E, x})^{\textup{ss}}$ sa semisimplification au sens de la décomposition de Jordan. Puisque~$\cL$ est $\iota$-pur de poids zéro, chaque élément $g$ dans $\rho(\Frob_{E, x})^{\textup{ss}}$ est semisimple avec des valeurs propres \textit{unitaires}, et l'adhérence de~$\langle g \rangle$ est donc un sous-groupe compact de $G_{\geom, \hspace{.3mm} \cL}(\CC)$. Il s'ensuit que $g$ appartient à un sous-groupe compact maximal et, comme ils sont tous conjugués, qu'il existe des éléments $h$ dans $K$ et $s$ dans~$G_{\geom, \hspace{.3mm} \cL}(\CC)$ tels que~$g=s^{-1} h s$. On affirme que $h$ est bien défini à conjugaison près par un élément de $K$. En effet, par le théorème de Peter--Weyl, les traces des représentations irréductibles de dimension finie~$\Lambda_K$ de $K$ séparent les classes de conjugaison. Or $\Lambda_K$ correspond d'après le théorème~\ref{theo:equiv} à une unique représentation $\Lambda$ du~$\alg_\ell$-groupe algébrique $G_{\geom, \hspace{.3mm} \cL}$~et, quel que soit le choix de~$h$, on a
\begin{displaymath}
\tr(\Lambda_K(h))=\tr(\Lambda(\rho(\Frob_{E, x})^{\textup{ss}}))=\tr(\Lambda(\rho(\Frob_{E, x})).
\end{displaymath} On notera $\theta_{E, x}$ la classe de conjugaison dans $K$ ainsi obtenue. 

Finalement, soit $\pi \colon K \to K^\#$ la projection canonique. Muni de la topologie quotient, l'espace $K^\#$ est compact et l'application $f \mapsto f_{\textup{cent}}=f \circ \pi$ permet d'identifier $\mathscr{C}(K^\#)$ aux fonctions \textit{centrales} continues sur $K$. L'image directe par $\pi$ de la mesure de Haar normalisée $\mu_K$ est une mesure de probabilité $\mu_{K^\#}$ sur $K^{\#}$ telle que
\begin{displaymath}
\int_{K^\#} f\mu_{K^\#}=\int_{K} f_{\textup{cent}}\hspace{.3mm} \mu_K.  
\end{displaymath} 

\begin{theo}[Deligne]\label{theo-deligne} Soit $\cL$ un faisceau $\ell$-adique lisse sur $X$. Supposons que~$\cL$ est $\iota$-pur de poids zéro et que les groupes $G_{\arith, \hspace{.3mm} \cL}$ et $G_{\geom, \hspace{.3mm} \cL}$ coïncident. Lorsque le degré des extensions $E/k$ tend vers l'infini, les classes de conjugaison $\{\theta_{E, x}\}_{x \in X(E)}$ s'équirépartissent dans $K^\#$ selon la mesure $\mu_{K^\#}$. 
\end{theo}

\noindent{\sc Démonstration}  --- Soit $f \colon K \to \CC$ une fonction centrale continue. Il faut démontrer que, quand $E$ parcourt les extensions de $k$ de degré assez grand pour que l'ensemble~$X(E)$ soit non vide, on a l'égalité 
\begin{equation}\label{eqn:amontrer}
\int_K f \mu_K=\lim_{|E| \to \infty} \frac{1}{|X(E)|} \sum_{x \in X(E)} f(\theta_{E, x}). 
\end{equation} 

D'après le théorème de Peter--Weyl, $f$ est limite uniforme de combinaisons linéaires finies de traces de représentations irréductibles de dimension finie $\Lambda_K$ de $K$; il suffit de traiter le cas d'une telle représentation. Comme $\mu_K$ a masse totale un, la fonction constante $f=1$ vérifie \eqref{eqn:amontrer} même sans passer à la limite. On peut donc supposer que~$\Lambda_K$ n'est pas triviale, auquel cas l'intégrale du membre gauche s'annule:  
\begin{displaymath}
 \int_K \tr_{\Lambda_K} \mu_K=0. 
\end{displaymath} 

Démontrons que le membre droit de \eqref{eqn:amontrer} est également nul (critère d'équirépartition  de Weyl). Invoquant à nouveau le théorème \ref{theo:equiv}, $\Lambda_K$ correspond à une unique représentation irréductible non triviale~$\Lambda$ du groupe algébrique $G_{\geom, \hspace{.3mm} \cL}$ sur~$\alg_\ell$. Compte tenu de l'hypothèse $G_{\arith, \hspace{.3mm} \cL}=G_{\geom, \hspace{.3mm} \cL}$, on peut composer $\Lambda$ avec
$$
\pi_1(X, \bar \eta) \stackrel{\rho}{\longrightarrow} G_{\geom, \hspace{.3mm} \cL}(\alg_\ell)
$$ et penser à $\Lambda \circ \rho$ comme à un système local $\ell$-adique sur~$X$. Notons-le $\Lambda(\cL)$ et regardons ses traces de Frobenius dans $\CC$ via $\iota.$ On a alors $\tr_{\Lambda_K}(\theta_{E, x})=\tr(\Frob_{E, x} \hspace{.3mm} | \hspace{.3mm} \Lambda(\cL))$, et la formule des traces de Grothendieck donne l'égalité 
\begin{displaymath}
\sum_{x \in X(E)} \tr_{\Lambda_K}(\theta_{E, x})=\sum_{i=0}^{2d} (-1)^i \tr(F_E \hspace{.3mm} | \hspace{.3mm} H^i_c(X_{\bar k}, \Lambda(\cL))). 
\end{displaymath}

Comme le groupe $G_{\geom, \hspace{.3mm} \cL}$ est réductif, $\rho$ est une représentation fidèle et $\Lambda$ est une représentation irréductible, il existe des entiers $a, b \geq 0$ tels que $\Lambda(\cL)$ soit sous-objet du système local \hbox{$\cL^{\otimes a} \otimes (\cL^\vee)^{\otimes b}$.} On en déduit que $\Lambda(\cL)$ est encore $\iota$-pur de poids zéro. La cohomologie à support compact $H^i_c(X_{\bar k}, \Lambda(\cL))$ est donc $\iota$-mixte de poids $\leq i$ d'après le théorème principal de Weil II. Par ailleurs, on a l'annulation 
$$
H^{2d}_c(X_{\bar k}, \Lambda(\cL))=(\Lambda(\cL)_{\bar \eta})_{\pi_1(X_{\bar k}, \bar \eta)}=0 
$$ car $\Lambda(\cL)_{\bar \eta}$ est une représentation irréductible non triviale de $\pi_1(X_{\bar k}, \bar \eta)$. Combinant ces informations, il vient
\begin{align}\label{eqn:estim}
\left|\sum_{x \in X(E)} \tr_{\Lambda_K}(\theta_{E, x})\right| &\leq \sum_{i=0}^{2d-1} \dim_{\alg_\ell} H^i_c(X_{\bar k}, \Lambda(\cL)) \cdot |E|^{i/2} \nonumber \\
&\leq \left(\sum_{i=0}^{2d-1} \dim_{\alg_\ell} H^i_c(X_{\bar k}, \Lambda(\cL)) \right) |E|^{(2d-1)/2}.
\end{align}

Par les inégalités de Lang--Weil, on a $2|X(E)| \geq |E|^{d}$ dès que le degré de l'extension~$E$ est assez grand pour que $|E|^{1/2} >2\sum_{i \leq 2d-1} b^i_c(X)$ et, en particulier, l'ensemble~$X(E)$ n'est pas vide. Joint à ce qui précède, ceci fournit la majoration
\begin{displaymath}
\left|\frac{1}{|X(E)|}\sum_{x \in X(E)} \tr_{\Lambda_K}(\theta_{E, x})\right| \leq \frac{2}{\sqrt{|E|}}\sum_{i=0}^{2d-1} \dim_{\alg_\ell} H^i_c(X_{\bar k}, \Lambda(\cL)). 
\end{displaymath} La limite du membre gauche est donc zéro lorsque le degré de $E$ tend vers l'infini.
\qed

\begin{rema}\label{rem:katz-sarnak} Dans \cite[Thm.\,9.2.6\,(3)]{katz-sarnak}, Katz et Sarnak démontrent qu'il existe une constante $C(X_{\bar k}, \cL)$, ne dépendant que de la variété $X_{\bar k}$ et du système local $\cL$, telle~que  
\begin{equation}\label{eqn:estim-katz-sarnak}
\sum_{i=0}^{2d} \dim_{\alg_\ell} H^i_c(X_{\bar k}, \Lambda(\cL)) \leq \dim(\Lambda) C(X_{\bar k}, \cL)
\end{equation} pour toute $\alg_\ell$-représentation de dimension finie $\Lambda$ de $G_{\geom, \hspace{.3mm} \cL}$. Si $X$ est une courbe, on peut d'après \cite[(3.6.2.1)]{katz-kloos} prendre  
\begin{equation*}\label{eqn:constant}
C(X_{\bar k}, \cL)=2g-2+N+\sum_{i=1}^N r_i
\end{equation*} où $g$ est le genre de la compactification lisse de $X$, $N$ le nombre de $\bar{k}$-points $x_1, \ldots, x_N$ à l'infini et~$r_i$ le plus grand saut dans la décomposition \eqref{eqn:dec-swan} de la monodromie locale de~$\cL$ en $x_i$. 
\end{rema}

\begin{rema} Le théorème \ref{theo-deligne} est un énoncé d'équirépartition verticale. Supposons maintenant que l'on se donne un schéma lisse $X$ sur un ouvert $S$ de $\Spec(\ZZ[1/\ell])$ et un système local $\ell$-adique $\iota$-pur de poids zéro $\cL$ sur $X$ tel que, pour chaque premier~$p$ dans~$S$, les groupes de monodromie géométrique et arithmétique du système local~$\cL_p$ sur~$X_{\FF_p}$ obtenu par réduction modulo $p$ coïncident et qu'ils sont tous égaux à un même groupe réductif. Au vu de la borne \eqref{eqn:estim-katz-sarnak}, la question si les classes de conjugaison~$\{\theta_{\FF_p, x}\}_{x \in X(\FF_p)}$ s'équirépartissent horizontalement est intimement liée au difficile problème de trouver des estimées explicites de la constante~$C(X_{\overline \FF_p}, \cL_p)$, par exemple de la majorer indépendamment de $p$, cf. \cite[App.]{katz-bull} et~\cite[\S9.6]{katz-sarnak}. 
\end{rema}

\subsection{Retour aux sommes de Kloosterman}

On donne maintenant l'application du théorème d'équirépartition  de Deligne aux sommes de Kloosterman. En plus des références originales de Deligne \cite{sommes-trigo} et Katz \cite{katz-kloos}, le lecteur pourra consulter avec profit le survol \cite{Laumon2000} de Laumon. Soient $\psi \colon \FF_q \to \alg_\ell^\times$ un caractère additif non trivial et $\cL_\psi$ le faisceau d'Artin--Schreier correspondant. Pour un entier $n \geq 2$, notons $\sigma \colon \GG_m^n \to \IA^1$ et $\pi \colon \GG_m^n \to \GG_m$ les morphismes somme et produit des coordonnées, respectivement, et posons 
\begin{displaymath}
\mathcal{K}l_n=R^{n-1}\pi_! \sigma^\ast \cL_\psi. 
\end{displaymath} 

D'après Deligne, $\mathcal{K}l_n$ est un faisceau $\ell$-adique lisse de rang $n$ sur $\GG_{m, \FF_q}$ ayant pour trace de Frobenius en un point $a \in \FF_q^\times$ la somme de Kloosterman au signe près: 
\begin{displaymath}
\tr(\Frob_{\FF_q, a} \hspace{.3mm} | \hspace{.3mm} \mathcal{K}l_n)=(-1)^{n-1} \Kl_n(a, q).
\end{displaymath} De plus, la flèche d'oubli des supports $R^{n-1}\pi_! \sigma^\ast \cL_\psi \to R^{n-1}\pi_\ast \sigma^\ast \cL_\psi$ est un isomorphisme~\cite[Thm.\,7.8]{sommes-trigo}. Comme $\sigma^\ast \cL_\psi$ est $\iota$-pur de poids zéro, le théorème principal de Weil II entraîne que $\mathcal{K}l_n$ est $\iota$-pur de poids $n-1$, d'où la majoration qui nous a permis de démontrer l'équirépartition  des angles des sommes de Gauss dans l'introduction: 
\begin{displaymath}
|\Kl_n(a, q)| \leq n q^{\frac{n-1}{2}}. 
\end{displaymath} 
 
Soit $\cL_n=\mathcal{K}l_n(\frac{n-1}{2})$ le système local $\ell$-adique $\iota$-pur de poids zéro obtenu en demi-tordant $\mathcal{K}l_n$ à la Tate. On sait que le déterminant $\det \cL_n=\wedge^n \cL_n$ est le système local trivial de rang un sur $\GG_{m, \FF_q}$ et que, pour $n$ pair ou pour $p=2$ et $n$ impair, on dispose d'un accouplement parfait $\cL_n \otimes \cL_n \to \alg_\ell$ qui est alterné dans le premier cas et symétrique dans le second. Il s'ensuit que le groupe de monodromie arithmétique~$G_{\arith, \hspace{.3mm} \cL_n}$ est inclus dans $\SL(n)$. Si $n$ est pair, c'est un sous-groupe du groupe symplectique $\mathrm{Sp}(n)$ et, si~$p=2$ et $n$ est impair, du groupe spécial orthogonal $\mathrm{SO}(n)$. Dans \cite[Thm.\,11.1]{katz-kloos}, Katz démontre que la monodromie est aussi grande que possible\footnote{D'après un théorème récent de Perret-Gentil \cite{corentin}, améliorant des résultats de Gabber \cite[Ch.\,12]{katz-kloos}, c'est encore le cas pour les groupes de monodromie \textit{entière}, qui sont des sous-groupes de $\Sp(n)$ ou $\SL(n)$ sur la complétion de l'anneau des entiers d'un corps cyclotomique en une place au-dessus de $\ell$, pourvu que $\ell$ soit plus grand qu'une constante effective dependant de $n$.}  en calculant
\begin{displaymath}
G_{\geom, \hspace{.3mm} \cL_n}=\begin{cases} \mathrm{Sp}(n) & n \text{ pair}, \\ \mathrm{SL}(n) & n \text{ impair et } p \neq 2, \\ \mathrm{SO}(n) & n \text{ impair}\neq 7 \text{ et } p=2, \\ G_2 & n=7 \text{ et } p=2, \end{cases}
\end{displaymath} où $G_2 \subset \mathrm{SO}(7)$ est l'un des groupes algébriques exceptionnels, défini comme le stabilisateur d'une certaine forme trilinéaire dans la représentation standard de dimension $7$ de $\mathrm{SO}(7)$. Avec un petit argument additionnel pour traiter le dernier cas \cite[\S11.3]{katz-kloos}, on en déduit l'égalité des groupes $G_{\geom, \hspace{.3mm} \cL_n}=G_{\arith, \hspace{.3mm} \cL_n}$.
 
Supposons $p \neq 2$ et notons $G$ le groupe de monodromie de $\cL_n$, qui est donc égal à $\Sp(n)$ si $n$ est pair et à $\SL(n)$ si $n$ est impair. L'intersection $K=G(\CC) \cap U(n)$ avec le groupe unitaire est un sous-groupe compact maximal de $G(\CC)$. Pour le système local~$\cL_n$, le nombre $C(\GG_{m, \bar{k}},\cL_n)$ de la remarque \ref{rem:katz-sarnak} vaut $1/n$ car $\cL_n$ a ramification modérée en~$0$ et ramification sauvage en $\infty$ avec $\Sw_\infty(\cL_n)=1$ et tous les sauts égaux à $1/n$. Par conséquent, \eqref{eqn:estim} et~\eqref{eqn:estim-katz-sarnak} donnent la majoration
\begin{displaymath}
\left|\frac{1}{q-1} \sum_{a \in \FF_q^\times} \tr_{\Lambda_K}(\theta_{\FF_q, a}) \right| \leq \frac{\dim(\Lambda_K)}{n} \frac{\sqrt{q}}{q-1}
\end{displaymath} pour toute représentation continue irréductible non triviale $\Lambda_K$ de $K$. Le résultat suivant, qui généralise le théorème \ref{thm:equi-kloos} dans l'introduction, en découle: 

\begin{theo}[Katz] Lorsque $q$ tend vers l'infini parmi n'importe quelle suite de puissances de nombres premiers impairs, les sommes de Kloosterman normalisées
\begin{displaymath}
\{\Kl_n(a, q)\slash q^{\frac{n-1}{2}}\}_{a \in \FF_q^\times}
\end{displaymath} se répartissent comme les traces de matrices aléatoires dans $K$. 
\end{theo} 
 
Il résulte de ce théorème que, pour chaque entier $r \geq 1$, la limite des moments
\begin{displaymath}
\lim_{q \to \infty} \frac{1}{q^{\frac{r(n-1)}{2}}(q-1)} \sum_{a \in \FF_q^\times} \Kl_n(a, q)^r
\end{displaymath} est égale à la multiplicité de la représentation triviale dans la puissance tensorielle~$r$\nobreakdash-ième de la représentation standard de $G$ correspondant à $\cL_n$. Par exemple, cette multiplicité vaut $2$ si $n=2$ et~\hbox{$r=4$,} ce qui explique le terme dominant dans le calcul \eqref{eqn:4ememoment} du quatrième moment des sommes de Kloosterman en une variable.  
  
\subsection{Le cas du groupe additif}\label{sec:additif} 

Dans ce paragraphe, on explique comment le théorème \ref{theo-deligne} permet d'étudier la répartition des transformées de Fourier, relatives à des caractères variables, de la fonction trace d'un faisceau $\ell$-adique sur la droite affine; c'est le pendant du théorème d'équirépartition de Katz pour le groupe additif. 

Soit $\mathbb{A}^1_k=\Spec k[x]$ et soit $M$ un objet de la catégorie dérivée $D^b_c(\mathbb{A}^1_k, \alg_\ell)$. Pour chaque extension finie $E$ de $k$ et chaque caractère additif $\rho \colon E \to \alg_\ell^\times$, on pose
\begin{displaymath}
S(M, E, \rho)=\sum_{x \in E} \rho(x)  \tr(\Frob_{E, x}\hspace{.3mm} |\hspace{.3mm} M).
\end{displaymath} Comment ces sommes se répartissent-elles lorsque le degré de $E$ tend vers l'infini? 

Moyennant le choix d'un caractère non trivial $\psi \colon k \to \alg_\ell^\times$, tous les caractères additifs de~$E$ sont de la forme $x \mapsto \psi(\tr_{E \slash k}(xy))$ pour un seul $y \in E$, de sorte que l'on peut récrire les sommes ci-dessus comme
\begin{displaymath}
S(M, E, y)=\sum_{x \in E} \psi(\tr_{E \slash k}(xy))  \tr(\Frob_{E, x}\hspace{.3mm} |\hspace{.3mm} M)
\end{displaymath} et les regarder comme étant paramétrées par la droite affine duale \hbox{$\mathbb{A}^1_k=\Spec k[y]$.} 

L'avantage de ce point de vue est que $y \mapsto -S(M, E, y)$ devient alors la fonction trace de la \textit{transformation de Fourier} 
\begin{align*}
\FT_\psi \colon D^b_c(\mathbb{A}^1_k, \alg_\ell) &\longrightarrow D^b_c(\mathbb{A}^1_k, \alg_\ell) \\ M &\longmapsto R\mathrm{pr}_{2,!}(\mathrm{pr}_1^\ast M \otimes \cL_{\psi(xy)}[1]), 
\end{align*}
où $\pr_1$ et $\pr_2$ désignent les projections de $\mathbb{A}^2_k=\Spec k[x, y]$ sur les deux facteurs. Introduite par Deligne dans les années 70, cette construction a été systematiquement étudiée par Brylinski, Katz et Laumon \cite{brylinski}, \cite{illusie}, \cite{katz-laumon}, \cite{laumon}, \cite{laumon-icm}. À l'instar de la transformation de Fourier classique, le foncteur $\FT_\psi$ est essentiellement involutif 
\begin{equation}\label{eqn:involutiviteFourier}
\FT_{\psi^{-1}} \circ \FT_{\psi}=(-1)
\end{equation} et échange le produit de convolution additive $M \star N=Rs_!(M \boxtimes N)$, où $s \colon \mathbb{A}^2_k \to \mathbb{A}^1_k$ désigne l'application somme, et le produit tensoriel:  
\begin{equation}\label{eqn:additiveconvolution}
\FT_\psi(M \star N)=\FT_\psi(M) \otimes \FT_\psi(N)[-1].
\end{equation}

On a défini la transformation de Fourier à l'aide de l'image directe à support compact car c'est celle-ci qui s'interprète en termes de sommes exponentielles par le biais de la formule des traces, mais on aurait pu également utiliser l'image directe usuelle. Le \og miracle de la transformation de Fourier \fg{} \cite[Thm.\,2.1.3]{katz-laumon} est qu'il n'y en a qu'une seule, c'est-à-dire que la flèche d'oubli des supports 
\begin{displaymath}
R\mathrm{pr}_{2,!}(\mathrm{pr}_1^\ast M \otimes \cL_{\psi(xy)}[1]) \longrightarrow R\mathrm{pr}_{2,\ast} (\mathrm{pr}_1^\ast M \otimes \cL_{\psi(xy)}[1])
\end{displaymath} est un isomorphisme. Voici quelques propriétés qui en résultent aussitôt: 

\begin{enumerate}
\item[(a)] Le foncteur $\FT_\psi$ commute à la dualité: $\DD \circ \FT_\psi=\FT_{\psi^{-1}} \circ \DD(1).$ 
\item[(b)]  Le foncteur $\FT_\psi$ induit une équivalence de catégories sur le faisceaux pervers:
$$\FT_\psi \colon \Perv(\mathbb{A}^1_k, \alg_\ell) \stackrel{\sim}{\longrightarrow} \Perv(\mathbb{A}^1_k, \alg_\ell).$$
 
\item[(c)] Le foncteur $\FT_\psi$ augmente de $1$ le poids des objets $\iota$-purs. 
\end{enumerate}

\noindent{\sc Démonstration}  --- La propriété (a) découle du fait que la dualité échange $R\mathrm{pr}_{2,!}$ et $R\mathrm{pr}_{2, \ast}$, ainsi que $\pr_1^\ast$ et~\hbox{$\pr_1^!=\pr_1^\ast(1)[2]$.} Grâce à l'involutivité \eqref{eqn:involutiviteFourier}, pour établir (b) il suffit de démontrer que la transformation de Fourier préserve les faisceaux pervers. Si~$M$ est pervers, alors~$\pr_1^\ast M[1] \otimes \cL_{\psi(xy)}$ est un faisceau pervers sur $\mathbb{A}^2_k$, et $\FT_\psi(M)$, vu comme image directe usuelle, est semipervers car le morphisme $\pr_2$ est affine; puisque son dual est égal à $\FT_{\psi^{-1}}(\DD(M)(1))$, il est aussi semipervers. Enfin, si $M$ est $\iota$-pur de poids $w$, alors~$\pr_1^\ast M \otimes \cL_{\psi(xy)}[1]$ est $\iota$-pur de poids $w+1$ et, comme les images directe usuelle et à support compact par la projection $\pr_2$ son isomorphes, (c) résulte du théorème principal de Weil II. \qed

La propriété (b) implique en particulier que la transformation de Fourier envoie des objets simples vers des objets simples. De plus, si $M$ n'est pas géométriquement isomorphe à un faisceau d'Artin--Schreier décalé $\cL_{\psi(ax)}[1]$, sa transformée de Fourier est de la forme $(j_\ast \cL)[1]$ pour un système local $\ell$-adique sur un ouvert $j \colon U \hookrightarrow \mathbb{A}^1_k$ et 
\begin{displaymath}
S(M, E, y)=\tr(\Frob_{E, y} \hspace{.3mm} | \hspace{.3mm} \cL)
\end{displaymath} pour tout $y \in U(E)$. Compte tenu de tout ce qui précède, le théorème d'équirépartition de Deligne se traduit en l'énoncé suivant:

\begin{theo} Soit $M$ un faisceau pervers géométriquement simple et $\iota$-pur de poids zéro sur \hbox{$\mathbb{A}^1_k$} qui n'est pas de la forme $\cL_{\psi(ax)}[1]$ et soit $(j_\ast \cL)[1]$ sa transformée de Fourier, où $\cL$ est un système local $\ell$-adique sur un ouvert $U$ de la droite affine duale. Supposons que les groupes de monodromie arithmétique et géométrique de $\cL$ coïncident et choisissons un sous-groupe compact maximal $K$ de ce groupe. Alors, les sommes~$\{S(M, E, y)\}_{y \in U(E)}$ se répartissent comme les traces de matrices aléatoires dans~$K$ lorsque le degré des extensions~$E \slash k$ tend vers l'infini. 
\end{theo}

Si l'on essaye d'adapter ces techniques à l'étude de la répartition des transformées de Fourier des fonctions trace d'un faisceau pervers sur le groupe multiplicatif $\mathbb{G}_{m, k}$, on est tout de suite confronté à l'obstacle qu'il n'y a pas de variété algébrique sur~$k$ dont les~$E$\nobreakdash-points paramètrent les caractères multiplicatifs $\chi \colon E^\times \to \alg_\ell$, ni \textit{a fortiori} de système local $\ell$-adique incarnant la transformation de Fourier relative à ces caractères\footnote{Par contre, le groupe des caractères $\ell$-adiques continus de $\pi_1^{\textup{mod}}(\GG_{m, \bar k})$ est muni d'une structure de $\alg_\ell$-schéma à une infinité dénombrable de composantes connexes isomorphes à $\Spec (\alg_\ell \otimes_{\ZZ_\ell}\ZZ_\ell[[X]])$, et l'on dispose d'un foncteur de $D^b_c(\GG_{m, \bar k}, \alg_\ell)$ vers la catégorie des modules à cohomologie bornée et cohérente sur ce schéma envoyant $M$ vers un objet dont la fibre en un caractère $\chi$ est $\mathrm{R}\Gamma(\mathbb{G}_{m, \bar k}, M \otimes \cL_\chi)$. C'est ce point de vue qui permet à Gabber et Loeser de démontrer certains des résultats du numéro suivant.}. Par ailleurs, on a vu avec les sommes de Gauss que leur répartition peut être gouvernée par des groupes comme $\GL(1)$ qui, contrairement aux groupes de monodromie géométrique d'un système local, ne sont pas semisimples. Katz a eu la belle idée de remplacer les points sur une variété algébrique par des foncteurs fibre sur une catégorie tannakienne comme espace de paramètres. Se rappelant que les systèmes locaux forment une catégorie tannakienne dans laquelle le groupe associé à un objet n'est rien d'autre que son groupe de monodromie, l'égalité \eqref{eqn:additiveconvolution} suggère de considérer à la place la catégorie des faisceaux pervers sur~$\GG_{m, k}$ munie du produit de convolution multiplicative. 

\section{Le théorème d'équirépartition  de Katz}\label{sec:multplicatif}

Ce numéro est consacré aux résultats de la monographie \cite{katz}. On introduit d'abord une catégorie tannakienne de faisceaux pervers sur le groupe multiplicatif~$\GG_m$ suivant Gabber et Loeser \cite{gabber-loeser}, ce qui permettra d'associer à ces objets des groupes algébriques sur~$\alg_\ell$ jouant le rôle des groupes de monodromie. Sous les bonnes hypothèses, il s'agit de groupes réductifs et l'on peut définir des classes de conjugaison~$\theta_{E, \chi}$ dans un sous\nobreakdash-groupe compact maximal $K$ de leurs points complexes. Le théorème de Katz affirme que les~$\theta_{E, \chi}$ s'équirépartissent selon  la mesure induite par la mesure de Haar normalisée. Après l'avoir démontré, on donnera quelques exemples de calcul du groupe tannakien.

\subsection{Bons caractères} 

Dans ce paragraphe et le suivant, on se place sur un corps algébriquement clos de caractéristique $p$, que l'on omettra de la notation. Soit $j \colon \GG_m \hookrightarrow \PP^1$ l'inclusion. 

\begin{defi}\label{defi:boncaracter} Soient $M$ un faisceau pervers sur $\GG_m$ et $\cL_\chi$ le faisceau de Kummer associé à un caractère continu $\chi$ du groupe fondamental modéré $\pi_1^{\textup{mod}}(\GG_m)$. On dit que~$\chi$ est {\rm bon} pour $M$ si la flèche canonique d'oubli des supports 
\begin{displaymath}
Rj_!(M \otimes \cL_\chi) \longrightarrow Rj_\ast(M \otimes \cL_\chi)
\end{displaymath} est un isomorphisme dans $D^b_c(\PP^1, \alg_\ell)$. Autrement, on dira que $\chi$ est {\rm mauvais} pour $M$. 
\end{defi}

\begin{lemm}\label{lem:nombrefini} Pour un objet donné $M$, il y a au plus $2\rg(M)$ mauvais caractères, où~$\rg(M)$ désigne le rang générique du faisceau constructible $\cH^{-1}(M)$. 
\end{lemm}

\noindent{\sc Démonstration}  --- Si $M$ est ponctuel, tous les caractères sont bons. En dévissant~$M$ en sa partie ponctuelle et non ponctuelle comme dans \eqref{eqn:devis}, on peut donc supposer que~$M$ est de la forme $\cF[1]$ pour un faisceau constructible $\cF$ sur $\GG_m$ qui n'a pas de sections ponctuelles. Posons $\cG=\cF \otimes \cL_\chi$. Comme le foncteur $j_!$ est exact, le caractère $\chi$ est bon pour $\cF$ si et seulement si $j_! \cG \cong j_\ast \cG$ et $R^1j_\ast \cG=0.$ Compte tenu de la description des foncteurs $j_!,$ $j_\ast$ et $R^1j_\ast$ dans l'exemple \ref{exemp:imagedirecte}, ces isomorphismes reviennent à dire que~$\cG$, en tant que représentation des groupes d'inertie $I_{\bar 0}$ et $I_{\overline \infty}$, n'a pas d'invariants ni de coinvariants. Or, si $j_0 \colon \GG_m \hookrightarrow \mathbb{A}^1$ désigne l'inclusion, les suites exactes
\begin{align}
&0 \longrightarrow \cG^{I_{\bar 0}} \longrightarrow H^0(\mathbb{A}^1, j_{0!} \cG[1]) \longrightarrow H^0(\GG_m, \cG[1]) \longrightarrow \cG_{I_{\bar 0}} \longrightarrow 0, \label{eqn:exactsequenceinvariants1} \\
&0 \longrightarrow \cG^{I_{\overline \infty}} \longrightarrow H^0_c(\GG_m, \cG[1]) \longrightarrow H^0(\mathbb{A}^1, j_{0!} \cG[1]) \longrightarrow \cG_{I_{\overline \infty}} \longrightarrow 0\label{eqn:exactsequenceinvariants2}
\end{align} montrent que ces invariants et coinvariants ont même dimension \cite[p.\,22]{katz}. Il suffit donc d'éviter les inverses des caractères modérés $\chi$ apparaissant dans la monodromie locale de $\cF$ autour de $0$ et~$\infty$ et il y en a au plus $2\rg(M)$. 
\qed

\subsection{Une catégorie tannakienne de faisceaux pervers sur $\GG_m$}

Soit $m \colon \GG_m \times \GG_m \to \GG_m$ la loi de multiplication. La catégorie dérivée $D^b_c(\GG_{m}, \alg_\ell)$ est munie de deux \textit{produits de convolution} 
\begin{displaymath}
M \star_! N=Rm_!(M \boxtimes N), \qquad  M \star N=Rm_\ast(M \boxtimes N), 
\end{displaymath} qui sont échangés par la dualité de Verdier. En combinant la suite spectrale de Leray pour le morphisme $m$, l'identité $m^\ast \cL_\chi=\cL_\chi \otimes \cL_\chi$ et la formule de K\"unneth \cite[\S2.5]{rigid}, on obtient des isomorphismes 
\begin{align}
\mathrm{R}\Gamma_c(\GG_m, (M \star_! N) \otimes \cL_\chi)&=\mathrm{R}\Gamma_c(\GG_m, M \otimes \cL_\chi) \otimes \mathrm{R}\Gamma_c(\GG_m, N \otimes \cL_\chi),  \label{eqn:Kunnethasupport} \\
\mathrm{R}\Gamma(\GG_m, (M \star N) \otimes \cL_\chi)&=\mathrm{R}\Gamma(\GG_m, M \otimes \cL_\chi) \otimes \mathrm{R}\Gamma(\GG_m, N \otimes \cL_\chi). \label{eqn:Kunnethsanssupport}
\end{align}

En général, aucun des produits de convolution ne préserve les objets pervers. Par exemple, le carré de convolution d'un objet de Kummer $\cL_\chi[1]$ est égal à $\cL_{\chi}[2]$ qui n'est pas pervers; dans un sens que l'on précisera ci-dessous, c'est la seule obstruction. Si $M$ et $N$ sont pervers, il en va de même pour $M \boxtimes N$ et, le morphisme $m$ étant affine, l'objet~$M \star N$ est semipervers. Comme $\DD(M \star_! N)=\DD(M) \star \DD(N)$ est également semipervers, si la flèche d'oubli des supports $M \star_! N \to M \star N$ est un isomorphisme, \textit{le} produit de convolution de $M$ et $N$ est pervers. Nous allons introduire une localisation de la catégorie $\Perv(\GG_m, \alg_\ell)$ sur laquelle cette flèche devient toujours un isomorphisme.

\begin{lemm}\label{lemm:euler} Soit $M$ un faisceau pervers sur $\GG_m$. 
\begin{enumerate}
\item[\rm (a)] On a $\chi(\GG_m, M) \geq 0.$
\item[\rm (b)] Si $\chi(\GG_m, M)=0$, alors $M$ est de la forme $\cF[1]$ pour un faisceau lisse $\cF$ qui est extension successive de faisceaux de Kummer $\cL_{\chi}$. 
\end{enumerate}
\end{lemm}

\noindent{\sc Démonstration}  --- La caractéristique d'Euler étant additive par rapport aux suites exactes, le dévissage \eqref{eqn:devis} de $M$ en sa partie ponctuelle et non ponctuelle donne 
\begin{displaymath}
\chi(\GG_{m}, M)=\chi(\GG_{m}, M_{\textup{pct}})+\chi(\GG_{m}, M_{\textup{npct}}).
\end{displaymath} Le premier terme est la dimension de $H^0(\GG_m, M_{\pct})$. Puisque $M_{\npct}=\cF[1]$ pour un faisceau constructible $\cF$ sur $\GG_{m}$ et que la caractéristique d'Euler de $\GG_{m}$ est nulle, la formule de Grothendieck--Ogg--Shafarevich est dans ce cas l'égalité
\begin{equation}\label{eqn:GOSGm}
\chi(\GG_{m}, M_{\textup{npct}})=\Sw_0(\cF)+\Sw_\infty(\cF)+\sum_{x \in \GG_m(\bar{k})} [\dr_x(\cF)+\Sw_x(\cF)]
\end{equation} dont tous les termes sont positifs car $\cF$ n'a pas de sections ponctuelles. 

Supposons maintenant $\chi(\GG_{m}, M)=0$. On a alors $H^0(\GG_{m}, M_{\textup{pct}})=0$, d'où $M_{\pct}=~0$ vu qu'un faisceau ponctuel sans sections globales est nul. Par suite, $M=\cF[1]$ et tous les termes du membre droit de \eqref{eqn:GOSGm} s'annulent. Puisque $\cF$ n'a pas de sections ponctuelles et que $\dr_x(\cF)=0$ pour tout $x$, il s'agit d'un faisceau lisse sur $\GG_m$. Comme $\Sw_0(\cF)$ et $\Sw_\infty(\cF)$ sont nuls, $\cF$ a ramification modérée à l'infini et correspond donc à une représentation continue de $\pi_1^{\textup{mod}}(\GG_m)$. Or ce groupe est abélien, comme on l'a vu dans le paragraphe~\ref{sec-ramification}, et ses caractères sont précisément les faisceaux de Kummer $\cL_\chi$. 
\qed

On dira qu'un faisceau pervers est \textit{négligeable} si sa caractéristique d'Euler est nulle. Soit $\text{\bf{Nég}}(\GG_m)$ la sous-catégorie pleine de $\Perv(\GG_m, \alg_\ell)$ formée des objets négligeables. Par le lemme \ref{lemm:euler}, elle est \textit{épaisse}: dans une suite exacte $0 \to A \to B \to C \to 0$, l'objet $B$ est négligeable si et seulement si $A$ et $C$ sont négligeables. On peut donc, d'après Gabber et Loeser, réaliser la construction suivante: soit $\overline{\text{Nég}}$ la sous-catégorie pleine de $D^b_c(\GG_m, \alg_\ell)$ formée des objets $M$ dont tous les faisceaux pervers de cohomologie~$\tensor[^p]{\cH}{^i}(M)$ sont négligeables. Alors $\overline{\text{Nég}}$ est une sous-catégorie épaisse, et la localisation  
\begin{equation}\label{eq-def-loc}
\Perv(\GG_m, \alg_\ell) \slash \text{\bf{Nég}}(\GG_m)
\end{equation} est équivalente à une sous-catégorie abélienne de la localisation de $D^b_c(\GG_m, \alg_\ell)$ relativement à $\overline{\text{Nég}}$, en fait le c\oe ur d'une $t$-structure \cite[\S3.6]{gabber-loeser}. 

\begin{prop} Soient $M$ et $N$ deux faisceaux pervers sur $\GG_m$. La flèche naturelle $M \star_! N \to M \star N$ dans $D^b(\GG_m, \alg_\ell)$ est un isomorphisme modulo \rm{$\overline{\text{Nég}}$}. 
\end{prop}

\noindent{\sc Démonstration}  --- Grâce au lemme \ref{lem:nombrefini}, il existe un caractère modéré $\chi$ tel que les complexes $\mathrm{R}\Gamma_c(\GG_m, \cdot)$ et $\mathrm{R}\Gamma(\GG_m, \cdot)$ soient isomorphes pour les trois objets $M \otimes \cL_\chi$, $N \otimes \cL_\chi$ et $(M \star_! N) \otimes \cL_\chi$. En particulier, la flèche naturelle
\begin{displaymath}
\mathrm{R}\Gamma_c(\GG_m, M \otimes \cL_\chi) \otimes \mathrm{R}\Gamma_c(\GG_m, N \otimes \cL_\chi) \to \mathrm{R}\Gamma(\GG_m, M \otimes \cL_\chi) \otimes \mathrm{R}\Gamma(\GG_m, N \otimes \cL_\chi)
\end{displaymath} est un isomorphisme. Or cette flèche s'identifie à 
\begin{displaymath}
\mathrm{R}\Gamma_c(\GG_m, (M \star_! N) \otimes \cL_\chi) \longrightarrow \mathrm{R}\Gamma(\GG_m, (M \star N) \otimes \cL_\chi)
\end{displaymath} d'après \eqref{eqn:Kunnethasupport} et \eqref{eqn:Kunnethsanssupport} et, comme la source est isomorphe à $\mathrm{R}\Gamma(\GG_m, (M \star_! N) \otimes \cL_\chi) $ par le choix de $\chi$, on en déduit que l'oubli des supports induit un isomorphisme
\begin{displaymath}
\mathrm{R}\Gamma(\GG_m, (M \star_! N) \otimes \cL_\chi) \longrightarrow \mathrm{R}\Gamma(\GG_m, (M \star N) \otimes \cL_\chi). 
\end{displaymath} 

Soit $C$ le cône du morphisme $M \star_! N \to M \star N$. Quitte à remplacer $C$ par $C \otimes \cL_\chi$, ce qui précède implique l'annulation $\mathrm{R}\Gamma(\GG_m, C)=0.$ Dans la suite spectrale
\begin{displaymath}
E_2^{a, b}=H^a(\GG_m, \tensor[^p]{\cH}{^b}(C)) \Longrightarrow H^{a+b}(\GG_m, C \otimes \cL_\chi), 
\end{displaymath} qui calcule la cohomologie d'un objet de la catégorie dérivée à partir de celle de ses faisceaux pervers de cohomologie, tous les termes $E_2^{a, b}$ avec $a \neq 0, 1$ s'annulent car~$\GG_m$ est une courbe affine. La suite spectrale dégénère donc à la deuxième page, et la condition~$\mathrm{R}\Gamma(\GG_m, C)=0$ équivaut à  $\mathrm{R}\Gamma(\GG_m, \tensor[^p]{\cH}{^b}(C))=0$ pour tout~$b.$ Par conséquent, les faisceaux pervers $\tensor[^p]{\cH}{^b}(C)$ sont négligeables. 
\qed

Il sera utile de disposer de représentants explicites des classes d'équivalence dans~\eqref{eq-def-loc}. 
 Katz considère la sous-catégorie $\mathbf{P}(\GG_m)$ de $\Perv(\GG_m, \alg_\ell)$ formée des faisceaux pervers n'ayant pas de sous-objets ni de quotients du type Kummer décalé $\cL_\chi[1]$. Dans~\hbox{\cite[\S2.6]{rigid},} il démontre que, pour un tel $M$, les objets $M \star_! N$ et $M \star N$ sont pervers quel que soit le faisceau pervers $N$ sur $\GG_m$. On peut alors parler de l'image de la flèche d'oubli des supports dans $\mathbf{P}(\GG_m)$ et définir la \textit{convolution intermédiaire}
\begin{displaymath}
M \star_{\textup{int}} N=\mathrm{im}(M \star_! N \longrightarrow M \star N).
\end{displaymath} En composant l'inclusion de $\bP(\GG_m)$ dans $\Perv(\GG_m, \alg_\ell)$ avec le passage à la localisation on obtient, d'après \cite[\S3.7]{gabber-loeser}, une équivalence de catégories 
\begin{displaymath}
\bP(\GG_m) \cong \Perv(\GG_m, \alg_\ell) \slash \text{\bf{Nég}}(\GG_m)
\end{displaymath} qui envoie la convolution intermédiaire $\star_{\textup{int}}$ sur \textit{le} produit de convolution. En particulier, chaque classe d'équivalence dans la localisation a un unique représentant dans $\bP(\GG_m)$. 
 
\skippointrait \begin{theo}[Gabber--Loeser, Katz, Deligne] 
\begin{enumerate}
\item[\rm (a)] Munie du produit de convolution intermédiaire~$\star_{\textup{int}}$, de l'objet neutre $\delta_1$ et de la dualité $M^\vee=[x \mapsto x^{-1}]^\ast \DD(M)$, la catégorie $\mathbf{P}(\GG_m)$ est tannakienne. 
\item[\rm (b)] Notons $j_0 \colon \GG_m \hookrightarrow \mathbb{A}^1_m$ l'inclusion. Pour n'importe quel faisceau de Kummer $\cL_\chi$ sur $\GG_m$, la correspondance
\begin{displaymath}
M \longmapsto H^0(\IA^1, j_{0!}(M \otimes \cL_\chi)) 
\end{displaymath} définit un foncteur fibre $\omega_\chi$ sur la catégorie tannakienne $\bP(\GG_m)$. 
\end{enumerate}
\end{theo}

\noindent{\sc Démonstration}  ---  La première assertion est le théorème 3.7.5 de Gabber--Loeser~\cite{gabber-loeser} et la deuxième, que Katz attribue à Deligne, est démontrée dans l'appendice à \cite{katz}. \qed

\begin{rema}\label{rema:dim-tan} Quel que soit le faisceau pervers $M$ sur $\GG_m$, la cohomologie~$H^i(\IA^1, j_{0!} M)$ s'annule en degré non nul. En effet, on peut supposer que $M$ est ou bien ponctuel, auquel cas l'assertion est évidente, ou bien extension intermediaire $\cF[1]$, auquel cas les groupes~$H^i(\IA^1, j_{0!} M)=H^{i+1}(\IA^1, j_{0!} \cF)$ s'annulent en degré~$i=1$ car $\IA^1$ est une courbe affine et en degré $i=-1$ car $\cF$ n'a pas de sections globales ponctuelles non nulles. Comme les caractéristiques d'Euler usuelle et à support compact d'un faisceau $\ell$-adique coïncident, on en déduit l'égalité~$\dim H^0(\IA^1, j_{0!}M)=\chi(\GG_m, M)$. En particulier, la dimension tannakienne d'un objet $M$ de $\bP(\GG_m)$ est donnée par sa caractéristique d'Euler.
\end{rema}

\subsection{Les classes de conjugaison $\theta_{E, \chi}$} 

Plaçons-nous maintenant sur un corps fini $k$. Soit $\bP_{\arith}(\GG_m)$ la sous-catégorie pleine de $\Perv(\GG_{m, k}, \alg_\ell)$ formée des objets $M$ tels que le faisceau pervers $M_{\bar k}$ sur $\GG_{m, \bar k}$ déduit de $M$ par extension des scalaires appartienne à la sous-catégorie $\bP(\GG_m)$ du paragraphe précédent. La convolution intermédiaire et le foncteur fibre 
$$
M \mapsto \omega(M)=H^0(\IA^1_{\bar k}, j_{0!}M)
$$ munissent $\bP_{\arith}(\GG_m)$ d'une structure de catégorie tannakienne. 

\begin{lemm}\label{lemme-weight0} La sous-catégorie pleine de $\bP_{\arith}(\GG_m)$ formée des faisceaux pervers~$\iota$\nobreakdash-purs de poids zéro est stable par le produit de convolution intermédiaire. 
\end{lemm}

\noindent{\sc Démonstration}  --- Si $M$ et $N$ sont des faisceaux $\iota$-purs de poids zéro sur $\GG_{m, k}$, leur produit extérieur~$M \boxtimes N$ est également $\iota$-pur de poids zéro sur $\GG_{m, k} \times \GG_{m, k}$. D'après le résultat principal de Weil II (théorème \ref{theo:del-cat}), $M \star_! N=Rm_!(M \boxtimes N)$ est $\iota$-mixte de poids~$\leq 0$. Puisque tout quotient d'un tel faisceau pervers est encore~$\iota$-mixte de poids~$\leq 0$ par~\hbox{\cite[Prop.\,5.3.1]{bbd},} on en déduit que $M \star_{\textup{int}} N$ est $\iota$-mixte de poids $\leq 0$. Or le dual~$\DD(M \star_{\textup{int}} N)$ s'identifie à la convolution intermédiaire $\DD(M) \star_{\textup{int}} \DD(N)$ et est donc~$\iota$\nobreakdash-mixte de poids~$\leq 0$ aussi. Par conséquent, $M \star_{\textup{int}} N$ est $\iota$-pur de poids zéro. 
\qed

Soit $M$ un faisceau pervers dans $\bP_{\arith}(\GG_m)$ semisimple et~$\iota$\nobreakdash-pur de poids zéro. En combinant les suites exactes \eqref{eqn:exactsequenceinvariants1} et \eqref{eqn:exactsequenceinvariants2} dans la preuve du lemme \ref{lem:nombrefini}, on voit que~\hbox{$Rj_! M \to Rj_\ast M$} est un isomorphisme si et seulement si les applications
\begin{displaymath}
H^0_c(\GG_{m, \bar k}, M) \longrightarrow \omega(M) \longrightarrow H^0(\GG_{m, \bar k}, M)
\end{displaymath} d'oubli des supports et de restriction sont des isomorphismes \cite[Thm.\,4.1]{katz}. Si c'est le cas, $\omega(M)$ est~$\iota$-pur de poids zéro d'après le théorème principal de Weil II, c'est-à-dire que les valeurs propres de Frobenius agissant sur $\omega(M)$ sont unitaires. 

Étant donnés une extension finie $E$ de $k$ et un caractère mutiplicatif $\chi \colon E^\times \to \alg_\ell^\times$, notons~$\cL_\chi$ le système local de Kummer sur $\GG_{m, E}$ défini dans l'exemple \ref{exem:lang}. Si $L$ est une extension finie de $E$, le tiré en arrière de $\cL_\chi$ sur $\GG_{m, L}$ est le faisceau de Kummer associé au caractère $\chi \circ N_{L \slash E}$. Vus comme des faisceaux lisses sur $\GG_{m, \bar k}$, les $\cL_\chi$ sont donc les caractères \textit{d'ordre fini} du groupe fondamental modéré 
\begin{displaymath}
\pi_1^{\textup{mod}}(\GG_{m, \bar k})=\varprojlim_{E/k \text{ finie}} E^\times , 
\end{displaymath} où les applications de transition sont données par la norme. On dira que le caractère $\chi$ est \textit{bon} ou \textit{mauvais} pour $M$ si $\chi$ est bon ou mauvais pour l'objet $M_{\bar k}$ sur $\GG_{m, \bar k}$ au sens de la définition \ref{defi:boncaracter}. Quitte à identifier $(E, \chi)$ avec $(L, \chi \circ N_{L / E})$ pour toute extension finie $L$ de $E$, le lemme~\ref{lem:nombrefini} entraîne qu'il y au plus $2 \rg(M)$ mauvais caractères pour $M$. 

En remplaçant $M$ par $M \otimes \cL_\chi$ dans ce qui précède, la caractérisation de bons caractères en termes du foncteur $\omega$ implique le résultat suivant:  

\begin{prop}\label{prop:fonctfib} Soient $M$ un faisceau pervers dans $\bP_{\arith}(\GG_m)$ et $\langle M \rangle^{\otimes}$ la sous\nobreakdash-catégorie tannakienne qu'il engendre. Supposons que $M$ est semisimple et $\iota$-pur de poids zéro. Si un caractère $\chi \colon E^\times \to \alg_\ell^\times$ est bon pour $M$, alors $\chi$ est bon pour n'importe quel objet $N$ dans $\langle M \rangle^{\otimes}$ et la correspondance
\begin{displaymath}
N \mapsto \omega_\chi(N)=H^0_c(\GG_{m, \bar k}, M \otimes \cL_\chi)
\end{displaymath} définit un foncteur fibre sur $\langle M \rangle^{\otimes}$. Le groupe $\omega_\chi(N)$ est $\iota$-pur de poids zéro et le Frobenius géométrique $F_E$ est un automorphsime tensoriel de $\omega_\chi$. 
\end{prop}

Fixons un foncteur fibre $\omega$ sur la catégorie tannakienne $\bP(\GG_m)$, par exemple
$$
\omega(M)=H^0(\IA^1_{\bar{k}}, j_{0!}M).
$$ Soient $M$ un faisceau pervers dans $\bP_{\arith}(\GG_m)$ et $M_{\bar{k}}$ l'objet correspondant de $\bP(\GG_m)$. Par le théorème principal des catégories tannakiennes,~$\omega$ induit des équivalences 
\begin{displaymath}
\langle M \rangle^{\otimes} \cong \Rep_{\alg_\ell}(\monarith) \qquad \langle M_{\bar{k}} \rangle^{\otimes} \cong \Rep_{\alg_\ell}(\mongeo)
\end{displaymath} entre la catégorie engendrée par $M$ (\textit{resp}. $M_{\bar{k}}$) et la catégorie des $\alg_\ell$-représentations de dimension finie du groupe des automorphismes tensoriels $\monarith$ (\textit{resp}. $\mongeo$) de la restriction de $\omega$ à $\langle M \rangle^{\otimes}$ (\textit{resp}. à $\langle M_{\bar k} \rangle^{\otimes}$). Il s'agit dans les deux cas de sous\nobreakdash-groupes de~$\GL(\omega(M))$. Concrètement, $\monarith$ est le sous-groupe des automorphismes linéaires~$\gamma$ tels que, pour toute construction tensorielle $M^{\ubar{r}, \ubar{s}}=\bigoplus M^{\otimes r_i} \otimes (M^\vee)^{\otimes s_i}$, l'action de~$\gamma$ sur~$\omega(M^{\ubar{r}, \ubar{s}})$ laisse globalement stables toutes les images des sous-quotients de~$M^{\ubar{r}, \ubar{s}}$. Comme l'objet $M_{\bar k}$ a plus de sous-quotients à respecter, on a l'inclusion
\begin{displaymath}
\mongeo \subseteq \monarith.
\end{displaymath}

Le corps de coefficients $\alg_\ell$ étant algébriquement clos, n'importe quel autre choix de foncteur fibre sur $\bP_{\arith}(\GG_m)$ donne lieu à un groupe isomorphe à $\monarith$ par un isomorphisme unique à un automorphisme intérieur près. Tout élément dans le groupe associé à un autre foncteur fibre définit donc une classe de conjugaison dans $\monarith$. En particulier, pour une extension finie $E$ de $k$ et un caractère multiplicatif $\chi \colon E^\times \to \alg_\ell^\times$, le Frobenius géométrique $F_E$ est un automorphisme tensoriel de $\omega_\chi$, d'où des classes de conjugaison dans le groupe $\monarith(\alg_\ell) \subset \monarith(\CC)$ que l'on notera $\Frob_{E, \chi}$. 

Si le faisceau pervers $M$ est $\iota$-pur de poids zéro, alors il est géométriquement semisimple d'après \cite[Thm.\,5.3.8]{bbd}. Comme la représentation fidèle $\mongeo \to \GL(\omega(M))$ correspondant à $M_{\bar k}$ est complètement réductible, $\mongeo$ est un groupe algébrique réductif. Dorénavant, on supposera que $M$ est \textit{arithmétiquement} semisimple, c'est\nobreakdash-à\nobreakdash-dire semisimple en tant qu'objet de $\bP_{\arith}(\GG_m)$. Ceci signifie que $M$ est un faisceau pervers semisimple sur $\GG_{m, k}$ qui, géométriquement, n'admet pas de faisceaux de Kummer décalés comme sous-objets ni comme quotients. Dans ce cas, le groupe $\monarith$ est réductif aussi. Si le caractère $\chi$ est bon pour $M$, la classe de conjugaison $\Frob_{E, \chi}$ a des valeurs propres unitaires d'après la proposition \ref{prop:fonctfib}. Choisissons un sous-groupe compact maximal $K$ du groupe de Lie réductif $\monarith(\CC)$. Reprenant mot par mot l'argument du paragraphe~\ref{sec:theodeligne}, on associe à chaque $\Frob_{E, \chi}$ une classe de conjugaison~$\theta_{E, \chi}$ dans~$K^{\#}$. 

\subsection{Le théorème principal}

\begin{theo}[Katz]\label{theo-prin} Soit $M$ un faisceau pervers dans $\bP_{\arith}(\GG_m)$. Supposons que~$M$ est $\iota$-pur de poids zéro, semisimple, et que les groupes $\mongeo$ et~$\monarith$ sont égaux. Choisissons un sous-groupe compact maximal $K$ de $\monarith(\CC)$. Lorsque le degré des extensions $E \slash k$ tend vers l'infini, les classes de conjugaison~$\{\theta_{E, \chi}\}_{\chi \hspace{.5mm} \text{bon}}$ s'équirépartissent dans $K^\#$ selon la mesure $\mu_{K^\#}$ induite par la mesure de Haar sur $K$. 
\end{theo}

\noindent{\sc Démonstration}  --- Soit $\Bon(E, M)$ l'ensemble de bons caractères $\chi \colon E^\times \to~ \alg_\ell^\times$ pour l'objet $M$; d'après le lemme \ref{lem:nombrefini}, il est non vide dès que le degré de $E$ est assez grand. Démontrons que, pour toute fonction centrale continue $f \colon K \to \CC$, on a
\begin{equation}\label{eq:proof-equidis}
\frac{1}{|\Bon(E, M)|}\sum_{\chi \in \Bon(E, M)} f(\theta_{E, \chi}) \longrightarrow \int_K f \mu_K  
\end{equation} lorsque le degré de $E$ tend vers l'infini. En utilisant le théorème de Peter--Weyl, on se réduit à prouver que le membre gauche de \eqref{eq:proof-equidis} converge vers zéro lorsque $f$ est la trace d'une représentation irréductible non triviale de dimension finie $\Lambda_K$ de $K$. 

Soit $\Lambda$ la seule $\alg_\ell$-représentation irréductible non triviale de $\monarith$ associée à~$\Lambda_K$ par le théorème \ref{theo:equiv}. Par le biais de l'équivalence $\langle M \rangle^\otimes \simeq \Rep_{\alg_\ell}(\monarith)$, elle correspond à un faisceau pervers \textit{simple} $N$ dans la sous-catégorie tannakienne de~$\bP_{\arith}(\GG_m)$ engendrée par $M$. De plus, l'hypothèse $\mongeo=\monarith$ entraîne que~$N$ est \textit{géométriquement simple}, et la non-trivialité de $\Lambda$ se traduit par le fait que~$N_{\bar{k}}$ n'est pas isomorphe à l'objet neutre $\delta_1$. En termes de $N$, la trace $\tr_{\Lambda_K}(\theta_{E, \chi})$ est égale à 
\begin{displaymath}
S(N, E, \chi)=\sum_{x \in E^\times} \chi(x) \tr(\Frob_{E, x}\hspace{.3mm} |\hspace{.3mm} N)
\end{displaymath} par la formule des traces de Grothendieck. Le résultat découlerait donc de l'estimée suivante lorsque le degré de $E$ tend vers l'infini:  
\begin{displaymath}
\frac{1}{|\Bon(E, M)|} \sum_{\chi \in \Bon(E, M)} S(N, E, \chi)=O\left(1/\sqrt{|E|}\right).   
\end{displaymath}

Pour la démontrer, remarquons d'abord qu'il suffit d'estimer la moyenne des sommes~$S(N, E, \chi)$ sur \textit{tous} les caractères multiplicatifs. En effet, en notant~$\Mauv(E, M)$ le complémentaire de $\Bon(E, M)$ dans $E^\times$, on a:  
\begin{align}\label{eqn:bon-mauvais}
\frac{1}{|\Bon(E, M)|} \sum_{\chi \in \Bon(E, M)} S(N, E, \chi)=\left(1+\frac{|\Mauv(E, M)|}{|\Bon(E, M)|} \right)\frac{1}{|E^\times|} \sum_{\chi} S(N, E, \chi) \nonumber \\ -\frac{1}{|\Bon(E, M)| }\sum_{\chi \in \Mauv(E, M)} S(N, E, \chi). 
\end{align} Or il y a au plus $2\rg(M)$ mauvais caractères pour l'objet $M$ et la majoration 
\begin{displaymath}
|S(N, E, \chi)| \leq \dim \omega(N)
\end{displaymath} est valable quel que soit $\chi$. En effet, comme $N$ appartient à la catégorie tannakienne engendrée par $M$, c'est un faisceau pervers semisimple et $\iota$-pur de poids zéro dans~$\bP_{\arith}(\GG_m)$; il en va de même pour $N \otimes \cL_\chi$. La cohomologie à support compact d'un tel faisceau est concentrée en degré zéro et $H^0_c(\GG_{m, \bar k}, N \otimes \cL_\chi)$ est $\iota$-mixte de poids~$\leq 0$, d'où l'inégalité voulue par la formule des traces et la remarque \ref{rema:dim-tan}. 

Estimons donc la moyenne des $S(N, E, \chi)$. Par orthogonalité des caractères,
\begin{displaymath}
\frac{1}{|E^\times|} \sum_{\chi} S(N, E, \chi)=\frac{1}{|E^\times|}\sum_{x \in E^\times} \tr(\Frob_{E, x} \hspace{.3mm} | \hspace{.3mm} N) \sum_{\chi} \chi(x)=\tr(\Frob_{E, 1} \hspace{.3mm} | \hspace{.3mm} N). 
\end{displaymath} D'après la classification des faisceaux pervers simples sur une courbe (exemple \ref{perv-courbes}), on sait que $N$ est de l'un de ces deux types: 
\begin{itemize}[wide =1em]
\item[$\bullet$] Si $N$ est ponctuel, alors il existe un point $t \in k^\times \setminus \{1\}$ et une unité $\ell$-adique $\alpha$ tels que $N=\alpha^{\deg} \otimes \delta_t$, la possibilité $t=1$ étant exclue par l'hypothèse que $N_{\bar{k}}$ n'est pas isomorphe à $\delta_1$. Dans ce cas, $\tr(\Frob_{E, 1} \hspace{.3mm} | \hspace{.3mm} N)=0$. 

\item[$\bullet$] Si $N$ n'est pas ponctuel, alors il est de la forme $\cF[1]$ pour un faisceau extension intermédiaire $\cF$ sur $\GG_m$. Comme $N$ est $\iota$-pur de poids zéro, $\cF$ est $\iota$-pur de poids $-1$, d'où les inégalités 
\begin{displaymath}
|\tr(\Frob_{E, 1} \hspace{.3mm} | \hspace{.3mm} N)|=|-\tr(\Frob_{E, 1} \hspace{.3mm} | \hspace{.3mm} \cF)| \leq \frac{\dim(\cF_{\bar{1}})}{\sqrt{|E|}} \leq \frac{\mathrm{rg}(\cF)}{\sqrt{|E|}}, 
\end{displaymath} où la dernière provient du fait que, le faisceau $\cF$ étant une extension intermédiaire, la fibre $\cF_{\bar 1}=(\cF_{\bar \eta})^{I_{\bar 1}}$ a dimension plus petite ou égale au rang générique de $\cF$. 
\end{itemize}

Dans les deux cas, on a $|\tr(\Frob_{E, 1} \hspace{.3mm} | \hspace{.3mm} N)| \leq \rg(N)\slash \sqrt{|E|}$. Compte tenu de \eqref{eqn:bon-mauvais}, on en déduit la majoration 
\begin{equation}\label{eqn:majorationpreuve}
\left|\frac{1}{|\Bon(E, M)|} \sum_{\chi \in \Bon(E, M)} S(N, E, \chi)\right| \leq \frac{2(\rg(N)+\dim \omega(N))}{\sqrt{|E|}}
\end{equation} pour toute extension $E\slash k$ de degré assez grand pour que $|\Mauv(E, M)| \leq \sqrt{|E|}-1$. C'est ce qu'il fallait pour conclure la démonstration. 
\qed

En prenant l'image directe de la mesure de Haar normalisée par la trace $\tr \colon K \to \CC,$ on trouve immédiatement le corollaire suivant: 

\begin{coro}\label{coro:katz} Sous les hypothèses du théorème \ref{theo-prin}, les sommes exponentielles~$\{S(M, E, \chi) \}_{\chi\,\text{bon}}$ se répartissent comme les traces de matrices aléatoires dans~$K$ lorsque le degré des extensions $E \slash k$ tend vers l'infini. 
\end{coro}

\begin{rema} Dans \cite[Thm.\,7.2]{katz}, Katz établit un résultat un peu plus général. Si $M$ est semisimple, $\mongeo$ est un sous-groupe \textit{distingué} de $\monarith$ et il suffit de supposer que le quotient $\monarith \slash \mongeo$ est isomorphe à un groupe cyclique $\ZZ/n\ZZ.$ Dans ce cas, on choisit un sous-groupe compact maximal $K_{\arith}$ de $\monarith(\CC)$ et l'on note~$K_{\geom}$ son intersection avec $\mongeo(\CC)$. C'est un sous-groupe distingué de $K_{\arith}$ avec quotient $K_{\arith} \slash K_{\geom}=\ZZ/n\ZZ$, d'où une application $K_{\arith}^\# \to \ZZ/n\ZZ$. Étant donné un entier $d$ modulo $n$, soient $K_{\arith, d}^\# \subset K_{\arith}^\#$ l'image inverse de $d$ par cette application et $\mu_{d}^\#$ la mesure de probabilité sur $K_{\arith, d}^\#$ déduite de la mesure de Haar normalisée sur~$K_{\geom}$. Le résultat est le suivant: pour n'importe quelle suite d'extensions~$E/k$ de degré congru à $d$ modulo $n$ dont les cardinaux tendent vers l'infini, les classes de conjugaison~$\{\theta_{E, \chi}\}_{\chi\,\text{bon}}$ s'équirépartissent dans $K_{\arith, d}^\#$ selon la mesure $\mu_{d}^\#$. 
\end{rema}

En présence de bornes uniformes pour les constantes intervenant dans l'estimée \eqref{eqn:majorationpreuve}, le théorème d'équirépartition de Katz admet la variante horizontale suivante: 

\begin{theo}[Katz]\label{thm:katzhorizontal} Soient $G$ un groupe réductif sur $\alg_\ell$ muni d'une représentation fidèle de dimension $n$ et $K$ un sous-groupe compact maximal de~$G(\CC)$. Soit $(k_i)$ une suite de corps finis de caractéristique distincte de $\ell$ dont les cardinaux tendent vers l'infini. Pour chaque $i$, soit $M_i$ un objet semisimple et $\iota$-pur de poids zéro de la catégorie~$\bP_{\arith}(\GG_{m, k_i})$ ayant dimension tannakienne~$n$. Supposons que l'égalité~$G_{\textup{geom}, M_i}=G_{\arith, M_i}=G$ soit satisfaite, l'objet $M_i$ correspondant à la représentation fidèle de $G$, et qu'il existe un nombre réel~$C \geq n$ tel que $\rg(M_i) \leq C$ pour tout $i$. Alors les classes de conjugaison $\{\theta_{k_i, \chi}\}_{\chi \in \Bon(k_i, M_i)}$ s'équirépartissent dans~$K^\#$. 
\end{theo}

\noindent {\sc Esquisse de démonstration} --- Soient $\Lambda$ une représentation irréductible non triviale de~$G$ et $N_i$ l'objet correspondant dans la sous-catégorie $\langle M_i \rangle^\otimes$ de $\bP_{\arith}(\GG_{m, k_i})$. D'après le critère d'équirépartition de Weyl, il suffit de démontrer que 
\begin{displaymath}
\frac{1}{|\Bon(k_i, M_i)|} \sum_{\chi \in \Bon(k_i, M_i)} S(N_i, k_i, \chi)
\end{displaymath} converge vers zéro lorsque $i$ tend vers l'infini. Comme il y a au plus $2\rg(M_i)$ mauvais caractères, \eqref{eqn:majorationpreuve} donne l'inégalité \begin{equation}\label{eqn:equirephorizontale}
\left|\frac{1}{|\Bon(k_i, M_i)|} \sum_{\chi \in \Bon(k_i, M_i)} S(N_i, k_i, \chi)\right| \leq \frac{2(\rg(N_i)+\dim \omega(N_i))}{\sqrt{|k_i|}}
\end{equation} pour tout $i$ tel que $2C \leq \sqrt{|k_i|}-1$. Le groupe $G \subset \GL(V)$ étant réductif, il existe des entiers $a$ et $b$ tels que $\Lambda$ soit sous-représentation de $V^{\otimes a} \otimes (V^\vee)^{\otimes b}$. Chaque faisceau pervers $N_i$ est donc sous-objet de la construction tensorielle $M_i^{\otimes a} \otimes (M_i^\vee)^{\otimes b}$ et
\begin{displaymath}
\dim \omega(N_i) \leq \dim(M_i^{\otimes a} \otimes (M_i^\vee)^{\otimes b})=n^{a+b} \leq C^{a+b}. 
\end{displaymath}
Par ailleurs, une analyse fine de la formule de Grothendieck--Ogg--Shafarevich montre que, pour tous faisceaux pervers semisimples $K$ et $L$ dans $\bP_{\arith}(\GG_{m, k_i})$, l'inégalité
\begin{displaymath}
\rg(K \star_{\textup{int}} L) \leq \rg(K) \dim \omega(L)+\rg(L) \dim \omega(K)
\end{displaymath} est satisfaite \cite[Thm.\,28.2]{katz}. Par récurrence sur $a+b$, on en déduit l'inégalité 
\begin{displaymath}
\rg(M_i^{\otimes a} \otimes (M_i^\vee)^{\otimes b}) \leq (a+b) (\dim \omega(M_i))^{a+b-1}\,\rg(M_i)
\end{displaymath} et donc $\rg(N_i) \leq (a+b) C^{a+b}$, ce qui permet de conclure que le membre gauche de \eqref{eqn:equirephorizontale} converge vers zéro lorsque $i$ tend vers l'infini. \qed

\subsection{Exemples}

Dans ce dernier paragraphe, nous expliquons comment déduire du théorème d'équirépartition  de Katz les résultats annoncés dans l'introduction. 

\begin{exem}[Sommes de Gauss] Soit $\psi \colon k \to \alg_\ell^\times$ un caractère additif non trivial. Alors $M=j_0^\ast \cL_\psi(1/2)[1]$ est un faisceau pervers simple sur $\GG_{m, k}$, pur de poids zéro, qui n'est pas géométriquement isomorphe à un faisceau de Kummer décalé (par exemple, car il a ramification sauvage à l'infini). Par la formule de Grothendieck--Ogg--Shafarevich, sa dimension tannakienne vaut 
\begin{displaymath}
\chi(\GG_m, M)=-\chi(\GG_m, j_0^\ast \cL_\psi)=\Sw_\infty(\cL_\psi)=1, 
\end{displaymath} d'où des inclusions $\mongeo \subseteq \monarith \subseteq \GL(1)$. Comme aucune puissance de convolution de $M_{\bar k}$ n'est l'objet neutre, le groupe $\mongeo$ est égal à $\GL(1)$ tout entier. Il s'ensuit que $M$ est un objet de $\bP_{\arith}(\GG_m)$ vérifiant les hypothèses du théorème. Soit~$E$ une extension de $k$ à $q$ éléments et posons $\psi_q=\psi \circ \tr_{E \slash k}$. Pour tout caractère multiplicatif~$\chi \colon E^\times \to \alg_\ell^\times$, on a 
\begin{displaymath} 
S(M, E, \chi)=-\frac{1}{\sqrt{q}} \sum_{x \in E^\times} \chi(x) \psi_q(x)=-\frac{g(\psi_q, \chi)}{\sqrt{q}}, 
\end{displaymath} de sorte que $\omega_\chi(M)$ est pur de poids zéro si $\chi$ n'est pas trivial et de poids $-1$ si $\chi$ est trivial. Les bons caractères pour $M$ sont donc les caractères non triviaux. Comme le compact maximal est $K=S^1$, le  corollaire \ref{coro:katz} est l'énoncé que les sommes de Gauss normalisées sont équiréparties dans $S^1$ selon la mesure de Haar.  
\end{exem}

Pour traiter les exemples des sommes d'Evans et de Rudnick, on aura besoin d'un résultat en théorie des représentations \cite[Thm.\,14.1]{katz}: 

\begin{prop}\label{prop-rep} Soit $M$ un faisceau pervers dans $\bP_{\arith}(\GG_m)$ de dimension tannakienne deux. Supposons que $M$ est pur de poids zéro, autodual et géométriquement semisimple. Si $M$ n'est géométriquement isomorphe à aucun de ses translatés non triviaux $[x \mapsto ax]^\ast M$, où $a \in \bar k \setminus \{1\}$, alors $\mongeo=\monarith=\SL(2).$
\end{prop}

\begin{exem}[Sommes d'Evans]\label{ex:evans} Soit $\psi \colon k \to \alg_\ell^\times$ un caractère additif non trivial. L'objet $M=\cL_{\psi(x-x^{-1})}(1/2)[1]$ est pur de poids zéro. Comme $\cL_{\psi(x-x^{-1})}$ est un faisceau de rang un à ramification sauvage à l'infini, $M$ est un objet géométriquement simple dans $\bP_{\arith}(\GG_m)$. Par la formule de Grothendieck--Ogg--Shafarevich, il a dimension
\begin{displaymath}
\chi(\GG_{m, \bar k}, M)=\Sw_0(\cL_{\psi(x-x^{-1})})+\Sw_\infty(\cL_{\psi(x-x^{-1})})=2. 
\end{displaymath} En écrivant $\cL_{\psi(x-x^{-1})}$ comme $\cL_{\psi(x)} \otimes \cL_{\psi(-1/x)}$, on voit que le dual de ce système local est $\cL_{\psi(x^{-1}-x)}$, d'où $M^\vee=M$ pour le dual tannakien. De plus, le translaté 
\begin{displaymath}
[x \mapsto ax]^\ast M=\cL_{\psi(ax-(ax)^{-1})}(1/2)[1]
\end{displaymath} n'est géométriquement isomorphe à $M$ pour aucun $a \neq 1$ (si tel était le cas, l'objet~$\cL_{\psi((a-1)x+(1-a^{-1})x^{-1})}(1/2)[1]$ serait géométriquement trivial, mais le système local sous-jacent a ramification sauvage à l'infini pour tout $a \neq 1$). Les hypothèses de la proposition \ref{prop-rep} sont donc satisfaites et~$\mongeo=\monarith=\SL(2)$. Comme $\cL_{\psi(x-x^{-1})}$ a ramification totalement sauvage en $0$ et $\infty$, il résulte de la preuve du lemme \ref{lem:nombrefini} que tous les caractères sont bons pour $M$. Le théorème d'équirépartition de Katz est dans ce cas l'énoncé que les angles~$\{\theta_{q, \chi}\}_{\chi}$ des sommes d'Evans s'équirépartissent dans l'intervalle~$[0, \pi]$ selon la mesure de Sato--Tate lorsque~$q$ tend vers l'infini. 
\end{exem}

\begin{exem}[Sommes de Rudnick]\label{ex:rudinick} Dans cet exemple, on suppose que la caractéristique de $k$ est impaire et on note $j_1 \colon \GG_{m, k} \setminus \{1\} \hookrightarrow \GG_{m, k}$ l'inclusion. Les sommes de Rudnick sont associées à l'objet
\begin{displaymath}
M=j_{1!} \cL_{\psi(\frac{x+1}{x-1})}(1/2)[1]
\end{displaymath} de la catégorie $\bP_{\arith}(\GG_m),$ qui est géométriquement simple et pur de poids zéro, par exemple parce que le faisceau $\cL_{\psi(\frac{x+1}{x-1})}$ a ramification totalement sauvage en $1$ et son prolongement par zéro coïncide donc avec l'image directe usuelle. La formule de Grothendieck--Ogg--Shafarevich donne dans ce cas l'égalité 
\begin{displaymath}
\chi(\GG_{m, \bar{k}}, M)=\dr_1(j_!\cL_{\psi(\frac{x+1}{x-1})})+\Sw_1(j_!\cL_{\psi(\frac{x+1}{x-1})}))=1+1=2
\end{displaymath} car $j_!\cL_{\psi(\frac{x+1}{x-1})}$ s'étend en un faisceau lisse en $0$ et $\infty$. De plus, l'identité \hbox{$\frac{x^{-1}+1}{x^{-1}-1}=-\frac{x+1}{x-1}$} implique que $M$ est autodual. Comme $1$ est le seul point de $\GG_m(\bar k)$ en lequel le faisceau~$j_!\cL_{\psi(\frac{x+1}{x-1})}$ n'est pas lisse, $M$ n'est géométriquement isomorphe à aucun translaté non trivial, d'où $\mongeo=\monarith=\SL(2)$ par la proposition \ref{prop-rep}. La monodromie de~$j_!\cL_{\psi(\frac{x+1}{x-1})}$ autour de $0$ et $\infty$ étant triviale, tout caractère \textit{non trivial}~$\chi$ est bon pour $M$ et le théorème d'équirépartition de Katz répond positivement à la question de Rudnick. 
\end{exem}

\bibliographystyle{smfplain}
    \bibliography{biblio}

\end{document}